\documentclass[12pt]{article}\usepackage{amsfonts,amssymb}

\usepackage{amsmath,bbm}
\usepackage{theorem}          
\usepackage{latexsym}
\usepackage{float}                   
\usepackage{times}
\usepackage{cite}

\newtheorem{tht}{Theorem}[section]

\newtheorem{thl}[tht]{Lemma}
\newtheorem{thp}[tht]{Proposition}

\newcommand{\anf}{\raisebox{0.2ex}{\scriptsize$\triangleleft$}}
\newcommand{\ang}{\raisebox{0.2ex}{\scriptsize$\triangleright$}}

\newcommand{\mn}{\medskip}    

\newcommand{\hsp}{\hspace{-1pt}} 
\newcommand{\hs}{\hspace{1pt}}

\newcommand{\lti}{\,{\scriptstyle\ltimes}\,} 
\newcommand{\rti}{\,{\scriptstyle\rtimes}\,} %
\newcommand{\sumop}{\mathop{\mbox{$\sum$}}}
\newcommand{\into}{\mathop{\mbox{$\int$}}}

\newcommand{\D}{{\mathcal{D}}} 
\newcommand{\E}{{\mathcal{E}}} 
 
\newcommand{\cN}{{\mathcal{N}}} 

\newcommand{\G}{{\mathcal{G}}} 
\newcommand{\Hh}{{\mathcal{H}}} 
\newcommand{\K}{{\mathcal{K}}}
\newcommand{\U}{{\mathcal{U}}}   
\newcommand{\T}{{\mathcal{T}}}

\newcommand{\X}{{\mathcal{X}}}
\newcommand{\A}{{\mathcal{A}}}

\newcommand{\cS}{\mathcal{S}} 
\newcommand{\cL}{\mathcal{L}}  
\newcommand{\V}{\mathcal{V}}

\newcommand{\Z}{\mathbb{Z}}
\newcommand{\N}{\mathbb{N}}
\newcommand{\R}{\mathbb{R}}
\newcommand{\C}{\mathbb{C}}

\newcommand{\Lin}{{\mathrm{Lin}}}

\newcommand{\Ker}{\mathrm{ker}\,}
\newcommand{\dd}{\mathrm{d}}

\newcommand{\rf}[1][]{\textup{\eqref{#1}}}

\newcommand{\SUU}{\cO(\mathrm{SU}_q(1,1))}

\newcommand{\suu}{\cU_q(\mathrm{su}_{1,1})}

\newcommand{\Eq}{\cO(\mathrm{E}_q(2))}
\newcommand{\ue}{\U_q(\mathrm{e}_2)}
\newcommand{\oc}{\cO(\C_q)}
\newcommand{\qd}{\cO(\mathrm{U}_q)}        
\newcommand{\rmU}{\mathrm{U}}

\newcommand{\cU}{\mathcal{U}}
\newcommand{\cO}{\mathcal{O}}
\newcommand{\cF}{\mathcal{F}}
\newcommand{\cM}{\mathcal{M}}

\newcommand{\dZ}{\mathbb{Z}}
\newcommand{\dR}{\mathbb{R}}
\newcommand{\dC}{\mathbb{C}}
\newcommand{\dN}{\mathbb{N}}
\newcommand{\dT}{\mathbb{T}}

\newcommand{\ip}[2]{\langle {#1},{#2}\rangle}
\title{Hilbert space representations of cross product algebras II}

\author{Konrad Schm\"udgen  and Elmar Wagner}

\begin{document}
\date{\small{Fakult\"at f\"ur Mathematik und 
Informatik\\ Universit\"at Leipzig, 
Augustusplatz 10, 04109 Leipzig, Germany\\ 
E-mail: 
Konrad.Schmuedgen@math.uni-leipzig.de /
Elmar.Wagner@math.uni-leipzig.de 
} }
\maketitle
\renewcommand{\theenumi}{\roman{enumi}}
\begin{abstract}\noindent
In this paper, we study and classify Hilbert space representations 
of cross product $\ast$-al\-ge\-bras of the quantized enveloping 
algebra $\ue$ with the coordinate algebras $\Eq$ of the 
quantum motion group and $\oc$ of the complex plane, 
and of the quantized enveloping algebra $\suu$ 
with the coordinate algebras $\SUU$ of the quantum group 
$\mathrm{SU}_q(1,1)$ and $\qd$ of the quantum disc. 
Invariant positive functionals and the corresponding Heisenberg 
representations are explicitely described.
\end{abstract}
{\small Keywords: Quantum groups, Unbounded representations}
%
%
%

%
\section{Introduction}
%
This is our  second paper on $\ast$-re\-pre\-sen\-ta\-tions of 
cross product $\ast$-al\-ge\-bras. While the first 
paper \cite{SW} deals mainly with coordinate algebras of 
compact quantum spaces, the present one is concerned with coordinate 
algebras of non-compact quantum spaces. We treat various cross 
product $\ast$-al\-ge\-bras related to the quantum motion group 
$\mathrm{E}_q(2)$ and to the quantum group $\mathrm{SU}_q(1,1)$. 
More precisely, we study cross product algebras of the 
Hopf $\ast$-al\-ge\-bra $\ue$ with the coordinate algebras $\oc$ of the 
quantum complex plane $\C_q$ and  $\Eq$ of the quantum group 
$\mathrm{E}_q(2)$, and of the Hopf $\ast$-al\-ge\-bra $\suu$ with 
the coordinate algebras $\qd$ of the quantum disc $\mathrm{U}_q$ 
and $\SUU$ of the quantum group $\mathrm{SU}_q(1,1)$.

The purpose of this paper is to study well-behaved Hilbert 
space representations of the cross product $\ast$-al\-ge\-bras, 
invariant positive functionals and  Heisenberg representations. 
Our main aim is to describe representations and invariant positive functionals by 
explicit formulas in terms of generators and functions of them. 

The quantum spaces considered in this paper are non-compact. 
Hence $\ast$-rep\-re\-sen\-ta\-tions of 
the their coordinate  algebras  involve unbounded operators. 
Thus, for $\ast$-re\-pre\-sen\-ta\-tions of the corresponding cross 
product algebras, elements of quantized enveloping algebras 
{\it and} of coordinate algebras act by unbounded operators. 
During the classification and derivation of representations,  
we occasionally add regularity conditions concerning 
the unbounded operators in order to exclude pathological behaviour. 
In other words, we classify only ``well-behaved'' 
$\ast$-re\-pre\-sen\-ta\-tions fulfilling these additional assumptions. 
In this paper, we use mainly the following 
regularity assumptions: For an algebraic relation $AB=BA$ 
with two hermitian elements $A$ and $B$ of a $\ast$-al\-ge\-bra, 
we assume that the corresponding Hilbert space operators are 
essentially self-adjoint and that 
their closures strongly commute (that is, their spectral projections 
mutually commute). If we have an algebraic relation $AB=pAB$ with $p$ 
real and $A$ hermitian, then we assume that $A$ is 
represented in the Hilbert space by an essentially 
self-adjoint operator and $\phi(\bar A)B\subset B\phi(p\bar A)$ 
for all $\phi\in\mathrm{L}^\infty(\R)$. If this is fulfilled, 
we say that the operator relation $AB=pBA$ holds in strong sense. 
As usual, $\bar A$ denotes the closure of an operator $A$. 
Finally, given a relation $N^\ast N=NN^\ast$ in a $\ast$-al\-ge\-bra, 
we require that $\bar N$ is normal when considered as a 
Hilbert space operator.

Our paper is organized as follows. In Section \ref{sec-pre}, 
we collect some preliminaries which are needed later. 
These are basic definitions and facts on cross product algebras, 
the definition of Heisenberg representations for non-unital 
$\ast$-al\-ge\-bras and two operator-theoretic lemmas. Representations 
of cross product algebras of $\ue$ and $\suu$ are investigated 
in Sections \ref{sec-Eq2} and \ref{sec-SU11}, respectively. 
After developing some useful algebraic properties of 
cross product algebras, all $\ast$-re\-pre\-sen\-ta\-tions satisfying our 
regularity assumptions are classified in terms of the actions of 
generators. At the end of both sections, invariant positive functionals for 
the corresponding quantum spaces and the Heisenberg representations are 
explicitely described. Since these quantum spaces are non-compact, 
invariant positive functionals are not finite on coordinate algebras.
Hence we extend the actions of $\ue$ and $\suu$ to larger 
function algebras $\cF(\X)$, where 
$\X=\mathrm{E}_q(2)$, $\C_q$, $\mathrm{SU}_q(1,1)$, $\mathrm{U}_q$, 
such that the invariant positive functionals are finite on appropriate 
subalgebras $\cF_0(\X)$. As a byproduct, we obtain explicit formulas 
for the actions of generators of $\ue$ and $\suu$ on general elements 
of coordinate and function algebras. 

For background information on these quantum spaces, we refer to 
\cite{SV,MMNNU,Wo4} for $\mathrm{SU}_q(1,1)$, 
\cite{KV,Wo3,Wo4,Koe} for  $\mathrm{E}_q(2)$ 
and \cite{KL,VK3} for the quantum disc. 
Invariant functionals on the corresponding quantum spaces 
appear in \cite{SV,KV,Koe,VK3}. 
A cross product algebra related to a differential calculus 
on $\oc$, its representations and an invariant functional 
for the quantum complex plane $\C_q$ were studied in \cite{CDP}.

In the reminder of this section, we set up some notation
and terminology. By a $\ast$-re\-pre\-sen\-ta\-tion 
of a $\ast$-al\-ge\-bra $\X$, 
we mean a homomorphism $\pi$ of $\X$ into the algebra of 
endomorphisms of a dense linear subspace $\D$ of a Hilbert space $\Hh$ 
such that 
$\langle\pi(x)\eta,\zeta\rangle=\langle\eta,\pi(x^\ast)\zeta\rangle$ 
for all $\eta,\zeta\in \D$ and $x\in\X$. Here 
$\langle\cdot,\cdot\rangle$ denotes the scalar product of $\Hh$. 
In order to simplify and shorten the notation, we shall drop the 
symbol $\pi$ and use the same letter, say $x$, for an element 
$x$ of the abstract $\ast$-al\-ge\-bra and the corresponding Hilbert space 
operator $\pi(x)$ under the representation $\pi$. 
This should cause no confusion. 
By another abuse of notation, we occasionally denote an operator and 
its closure by the same symbol. Throughout we shall describe 
$\ast$-re\-pre\-sen\-ta\-tions by formulas for the actions of algebra 
generators. In all cases, an invariant dense domain $\D$ is 
easily constructed (for instance, by taking the linear span of 
base vectors). 

For a self-adjoint operator $T$, 
the notation $\sigma (T) \sqsubseteq (a,b]$ 
means that the spectrum $\sigma(T)$ of $T$ is contained in $[a,b]$ and 
that $a$ is not an eigenvalue. 
The notations $\sigma (T) \sqsubseteq [a,b)$ and $T>0$ 
have a similar meaning.   
Let $I$ be an index set, $\K$ be a Hilbert space 
and $\Hh=\oplus_{\iota\in I} \Hh_\iota$, 
where each $\Hh_\iota$ is the same Hilbert space $\K$.   
With $\eta$ in $\K$ and $j$ an index, $\eta_{j}$ denotes 
the vector in $\Hh$  which has $\eta$
in the $j$th position and zero elsewhere.  
We set $\eta_{j}=0$ if $j\notin I$.

Unless stated otherwise, $q$ stands for a real number belonging to 
the open interval $(0,1)$. We abbreviate $\lambda:=q-q^{-1}$
and 
$$\lambda_n:=(1-q^{2n})^{1/2},\quad  
\alpha_k(A):=(1+q^{2k}A^2)^{1/2},\quad 
\beta_k(A):=(1+q^{-2k}A^{-2})^{1/2},
$$
where $n\in\N_0$, $k\in\Z$ and $A$ is a self-adjoint operator.
%
%
\section{Preliminaries}
                                                \label{sec-pre}
%
\subsection{Basics of cross product algebras}
             
                                                \label{basics}
This subsection reviews basic definitions and 
facts on cross product algebras. 
For further details, see \cite{SW}. We shall denote the comultiplication, 
the counit and the antipode of a Hopf algebra 
by $\Delta$, $\varepsilon$ and $S$, respectively, and use 
the Sweedler notations 
$\varphi(x)=x_{(1)}\otimes x_{(2)}$ for a coaction $\varphi$ and 
$\Delta(x)=x_{(1)}\otimes x_{(2)}$.

Let $\U$ be a Hopf $\ast$-al\-ge\-bra. 
A $\ast$-al\-ge\-bra $\X$ is called a {\it right $\cU$-mod\-ule 
$\ast$-al\-ge\-bra} 
if $\X$ is a right $\cU$-mod\-ule with action $\anf$ satisfying 
$$
   (xy)\anf f = (x\anf f_{(1)})(y\anf f_{(2)}),\quad 
   (x\anf f)^\ast = x^\ast \anf S(f)^\ast,\quad x,y\in X,\ f\in\cU. 
$$
For an algebra $\X$ with unit $1$, we additionally require 
$1\anf f=\varepsilon(f)1$ for $f\in\cU$.

The {\it right cross product algebra $\cU\lti\X$} is the 
linear space $\cU\otimes \X$ equipped with product and involution
defined by 
\begin{equation}                                      \label{cross0}
(g\otimes x)(f\otimes y)= gf_{(1)}\otimes (x\anf f_{(2)})y,\quad  
 (f\otimes x)^\ast=f^\ast_{(1)}\otimes (x^\ast\anf f^\ast_{(2)}),
\end{equation} 
where $x,y\in\X$ and $g,f\in\cU$.
Let $\U_0$ be a $\ast$-sub\-al\-ge\-bra of $\U$ which is a right coideal 
of the Hopf algebra $\U$ (that is, $\Delta(\U_0)\subseteq \U_0\otimes\U$). 
From (\ref{cross0}), it follows that the linear 
subspace $\U_0\otimes\X$ of 
$\U\otimes \X$ is a $\ast$-sub\-al\-ge\-bra of $\U\lti\X$. 
We shall denote this subalgebra  by $\U_0\lti\X$.

If $\X$ and $\U_0$ have unit elements, we can consider $\X$ and $\U_0$ 
as subalgebras of $\U_0\lti\X$ by identifying $f\otimes 1$ and $1\otimes x$ 
with $f$ and $x$, respectively. 
In this way,  $\U_0\lti\X$ can be viewed as the 
$\ast$-al\-ge\-bra generated by the two subalgebras $\U_0$ and $\X$ 
with respect 
to the cross commutation relations
\begin{equation}                                       \label{cross1}
xf=f_{(1)} (x\anf f_{(2)}),\quad x\in\X,\ f\in\U_0.
\end{equation}

Let $\A$ be a Hopf $\ast$-al\-ge\-bra, $\ip{\cdot}{\cdot}$ 
a dual pairing of Hopf $\ast$-al\-ge\-bras $\U$ and 
$\A$, and $\X$ a left $\A$-comodule $\ast$-al\-ge\-bra.
Then $\X$ becomes a right 
$\U$-mod\-ule $\ast$-al\-ge\-bra with right \mbox{action $\anf$} given by 
\begin{equation}                                        \label{rechts}
  x\anf f = \ip{f}{x_{(1)}}x_{(2)},\quad x\in\X,\ f\in\U.
\end{equation}
In this case, Equation \rf[cross1] reads  
\begin{equation}                                       \label{cross2}
xf=f_{(1)} \langle f_{(2)},x_{(1)}\rangle x_{(2)},\quad x\in\X,\ f\in\U.
\end{equation}

The above definitions have their left-handed 
counterparts. Suppose that $\X$ is a 
left $\U$-mod\-ule $\ast$-al\-ge\-bra, that is, $\X$ is a 
left $\U$-mod\-ule with action $\ang$ satisfying 
$$
   f\ang (xy)=(f_{(1)}\ang x)(f_{(2)}\ang y),\quad
   (f\ang x)^\ast  =S(f)^\ast \ang x^\ast ,\quad x,y\in X,\ f\in\cU,
$$
and $f\ang 1=\varepsilon(f)1$ for $f\in\cU$ if $\X$ has a unit $1$.
Then the vector space $\X\otimes \U$ is a $\ast$-al\-ge\-bra, called 
{\it left cross product algebra} and denoted by $\X\rti\U$, 
with product and involution defined by
$$
(y\otimes f)(x\otimes g) \hsp=\hsp y(f_{(1)}\ang x)\otimes f_{(2)}g,  \ \
(x\otimes f)^\ast \hsp=\hsp (f_{(1)}^\ast \ang x^\ast) \otimes f^\ast_{(2)},
\ \ \, x,y\in\X,\ f,g\in\U.
$$

When $\X$ has a unit, $\X\rti\U$ can be 
considered as the $\ast$-al\-ge\-bra generated by the subalgebras $\X$ and 
$\U$ with cross relations 
\begin{equation}                                       \label{oma}
fx=(f_{(1)}\ang x)f_{(2)},\quad x\in\X,\  f\in\U.
\end{equation}

If $\X$ is a right  comodule $\ast$-al\-ge\-bra of a 
Hopf $\ast$-al\-ge\-bra $\A$ and 
$\langle \cdot,\cdot\rangle$ is 
a dual pairing of Hopf $\ast$-al\-ge\-bras $\A$ and $\U$, 
then $\X$ is a left $\U$-mod\-ule 
$\ast$-al\-ge\-bra with left action $\ang$ given by 
\begin{equation}                                       \label{ang}
  f\ang x=x_{(1)}\ip{f}{x_{(2)}},\quad x\in\X,\ f\in\cU,
\end{equation}
and Equation \rf[oma] can be written 
\begin{equation}                                     \label{cross3}
fx = x_{(1)} \langle f_{(1)},x_{(2)} \rangle f_{(2)},
\quad x\in\X,\ f\in\U.
\end{equation}

In many cases, it suffices to consider the right-handed version only. 
The following simple lemma shows how one can pass 
under certain conditions 
from a left to a right action. 
\begin{thl}                                           \label{r-l}
Let $\U$ be a Hopf $\ast$-al\-ge\-bra and $\X$ a 
left\ \,$\U$-mod\-ule $\ast$-al\-ge\-bra with left action $\ang$.
Suppose $\phi :\U\rightarrow\U$ is an algebra 
anti-automorphism and a coalgebra homomorphism, that is, $\phi$ 
is a bijective linear map satisfying  
$\phi(fg)=\phi(g)\phi(f)$ and 
$\Delta\circ \phi(f)=(\phi\otimes\phi)\circ\Delta (f)$   
for all $f,g\in\U$. 
Assume that $\ast\circ S\circ\phi=\phi\circ \ast\circ S$. 
Then the formula 
\begin{equation}                                      \label{ract}
    x\anf f:= \phi(f)\ang x,\quad f\in\U,\ x\in\X,      
\end{equation}
defines a right $\U$-ac\-tion on $\X$ which turns $\X$ into a 
right $\U$-mod\-ule $\ast$-al\-ge\-bra.
\end{thl}

The idea of the next lemma is taken from the paper \cite{F}.

\begin{thl}                      \label{2.2}
Let $\U\lti\X$ be a right cross product algebra. 
Let $\V$ be a right coideal of\ \,$\U$ and let $\X_0$  be 
a set of generators of the algebra $\X$. 
Suppose that there exists a linear mapping 
$\rho:\V\rightarrow \X$ such that
\begin{equation}                      \label{xvcon}
x\rho(v)=\rho(v_{(1)})(x\anf v_{(2)})
\end{equation}
for $x\in\X_0$ and $v\in\V$. Then, for each $v\in\V$, 
the element $\xi(v) := \rho (v_{(1)}) S(v_{(2)})$ 
commutes with the algebra $\X$ in $\U\lti\X$.
\end{thl}

\noindent
{\bf Proof.} 
First observe  that if (\ref{xvcon}) holds 
for $x$ and $y$ in $\X$, then it holds also for the product $xy$. 
Thus, we can assume that $\X_0=\X$. Let $x\in\X$ and $v\in\V$. 
Using Equations (\ref{cross1}) and (\ref{xvcon}), we compute
\begin{align*}
\quad  x\xi (v)&=x\rho (v_{(1)}) S(v_{(2)})
=\rho(v_{(1)})(x\anf v_{(2)}) S(v_{(3)}) &\\ 
 &=\rho (v_{(1}) S(v_{(4)}) (x\anf v_{(2)} S(v_{(3)})
=\rho(v_{(1)}) S(v_{(2)}) (x\anf 1)  = \xi (v)x.& \qquad \Box
\end{align*}
\subsection{Heisenberg representations}
                                          \label{basics-Heis}
Let $\U$ be a Hopf $\ast$-al\-ge\-bra and let 
$\X_1$ and $\X_0$ be right $\U$-mod\-ule $\ast$-al\-ge\-bras such that $\X_0$ 
is a left $\X_1$-mod\-ule satisfying
\begin{equation}                                     \label{bimod}
 (a. x)^\ast y=x^\ast (a^\ast. y),\quad
(a. y)\anf f=( a\anf f_{(1)}). ( y\anf f_{(2)})
\end{equation}
for $x,y\in\X_0$, $a\in\X_1$ and $f\in\U$. Here, $a. x$ stands for the 
left action of $a\in\X_1$ on $x\in\X_0$. 
Note that the first condition of \rf[bimod] appeared 
in \cite{S-hep}. 
Suppose that $h$ is a positive 
linear functional on $\X_0$ (i.e., $h(x^\ast x) \ge 0$ for all $x\in\X_0$) 
which is $\U $-in\-var\-i\-ant (i.e., $h(x\anf f)= \varepsilon (f)h(x)$ 
for $x\in\X_0$ and $f\in\U )$. 
We associate with such a functional $h$ a unique 
$\ast$-re\-pre\-sen\-ta\-tion $\pi_h$ of 
$\U \lti\X_1$ called 
the {\it Heisenberg representation associated with} $h$. 
Its construction is similar to the one for 
unital $\ast$-al\-ge\-bras (see \cite{SW}). 
Let us review the basic ideas. 
Set $\cN=\{x\in\X_0\,;\,h(x^\ast x)=0\}$ and 
$\tilde\X_0=\X_0/\cN$. We write $x\mapsto \tilde x$ to denote 
the canonical mapping $\X_0\rightarrow \tilde\X_0$. 
Then the linear space $\tilde\X_0$ is an inner product space with inner 
product 
\begin{equation}                                            \label{ip}
\langle \tilde{x},\tilde{y}\rangle=h(y^\ast x), \quad x,y\in\X_0.
\end{equation}
The action of the cross product algebra $\U \lti\X_1$ on $\tilde\X_0$ 
is given by 
$$
\pi_h(f\otimes a)\tilde{x} := ((a. x) \anf S^{-1} (f))\tilde{}\hs\hs, 
\quad x\in\X_0,\   a\in\X_1,\  f\in\U .
$$
That $\pi_h$ is a well defined $\ast$-re\-pre\-sen\-ta\-tion of the 
$\ast$-al\-ge\-bra $\U \lti\X_1$ has been proved in 
\cite[Proposition 5.3]{SW} for unital $\ast$-al\-ge\-bras $\X_0=\X_1$. 
With some necessary modifications, the proof remains valid in the present  
situation as well. 
If $\X_0\subset\X_1$, then $\X_0$ is a $\ast$-sub\-al\-ge\-bra of 
$\U \lti\X_1$ and the restriction of $\pi_h$ to $\X_0$ is just the 
GNS-re\-pre\-sen\-ta\-tion associated with $h$. 
In the sequel, we write simply $x$ instead of $\tilde x$. 

Condition  \rf[bimod] is satisfied if there is a 
right $\U$-mod\-ule $\ast$-al\-ge\-bra $\X$ such that $\X_1$ and $\X_0$ are 
right $\U$-mod\-ule $\ast$-sub\-al\-ge\-bras and $\X_0$ is an ideal of $\X$. 
In all our examples below, $\X$ will be  a $\ast$-al\-ge\-bra 
of functions on a non-compact quantum space which contains the 
coordinate algebra, here denoted by $\X_1$, as a subalgebra and  
$\X_0$ can be considered as a $\ast$-al\-ge\-bra of functions 
with compact support. In this situation, the left action 
$a. x$ of $a\in\X_1$ on $x\in\X_0$ becomes the product $ax$
in the algebra $\X$, and the restriction of the Heisenberg representation 
of the cross product $\ast$-al\-ge\-bra $\U \lti\X$
to the $\ast$-sub\-al\-ge\-bra $\U \lti\X_1$ is just the 
Heisenberg representation of $\U \lti\X_1$ 
(both associated with the same functional $h$ on $\X_0$). 
%
%
\subsection{Two auxiliary lemmas}
%
%
\begin{thl}                                     \label{1}
Let $p\in(0,1)$ and $m\in\N_0$. 
Suppose that $w$ is a unitary operator and 
$A_0,\ldots,A_m$ are self-adjoint operators 
on a Hilbert space $\Hh$. 
Assume that 
\begin{equation}
 A_0>0,\quad wA_0 \subseteq pA_0w,      \label{wA}
\end{equation}
and, if $m>0$, 
\begin{equation}                         \label{wAk}
 wA_k \subseteq A_kw,\ \, k=1,\ldots, m,\quad 
A_0,\ldots,A_m\  \mbox{strongly\ commute}. 
\end{equation}
Then, up to unitary equivalence, there exist a Hilbert space 
$\Hh_0$, strongly commuting self-adjoint operators 
$B_0,\ldots,B_m$ on $\Hh_0$,  
and a dense subspace $\D_0\subseteq\Hh_0$
such that 
$\Hh=\oplus^\infty_{n=-\infty} \Hh_n$, 
where each $\Hh_n$ is $\Hh_0$, 
$\D:=\Lin\{\eta_n\,;\,n\in\Z,\,\eta\in\D_0\}$ is invariant under 
$w$, $A_0,\ldots,A_m$, the operator $B_0$ satisfies 
$\sigma (B_0) \sqsubseteq (p,1]$, and the actions of 
$w$, $A_0,\ldots,A_m$ on $\Hh$ are determined by 
$$
 w\eta_n=\eta_{n-1},\quad A_0\eta_n=p^nB_0\eta_n,\quad 
A_k\eta_n=B_k\eta_n,\ \,k=1,\ldots,m,\quad \eta\in\D_0.
$$
In particular, 
$$
   wA_0=pA_0w,\quad wA_k= A_kw,\ \,k=1,\ldots,m,\quad A_jA_i=A_iA_j\quad 
\mbox{on}\ \,\D.
$$
If $A_k$ satisfies $\sigma (A_k) \sqsubseteq (a,b]$ or 
$\sigma (A_k) \sqsubseteq [a,b)$, where $a,b\in\R\cup \{\pm\infty\}$, $a<b$,   
and $k\in\{1,\ldots,m\}$, then the same holds for $B_k$. 
\end{thl}
{\bf Proof.}
Let $e(\mu)$ denote the spectral projections of $A_0$.
By \rf[wA], the self-adjoint operators $wA_0w^\ast$ and $pA_0$
coincide, hence $we(\mu)w^\ast=e(p^{-1}\mu)$. 
Define $\Hh_n=\big(e(p^n)-e(p^{n+1})\big)\Hh$. 
Then $w\Hh_n= \big(e(p^{n-1})-e(p^{n})\big)w\Hh=\Hh_{n-1}$. 
After applying an obvious unitary transformation, we 
may assume that $\Hh_n=\Hh_0$ and $w\eta_n =\eta_{n-1}$ 
for $\eta \in \Hh_0$.
Denote by $B_0$ the restriction of $A_0$ to $\Hh_0$. 
Then, by the definition of $\Hh_0$, $\sigma (B_0) \sqsubseteq (p,1]$
and $A_0\eta_n =p^n w^nA_0w^{n\ast}\eta_n= p^nw^nA_0\eta_0 
=p^nB_0\eta_n$. 

Let $m>0$ and $k\in\{1,\ldots,m\}$. The operator $A_k$ commutes strongly 
with $A_0$ and thus it commutes with the spectral projections of $A_0$. 
Therefore $A_k$ leaves each Hilbert space $\Hh_n\,(=\Hh_0)$ invariant. 
Denote by $A_{kn}$ the restriction of $A_k$ to $\Hh_n$. 
From \rf[wAk], we conclude that $A_{kn}\eta_{n-1}=A_{k,n-1}\eta_{n-1}$ 
for all $\eta_n$ in the domain of $A_{kn}$. 
Since both $A_{kn}$ and $A_{k,n-1}$ 
are self-adjoint operators on $\Hh_0$, we have $A_{kn}=A_{k,n-1}$, 
so all $A_{kn}$ are equal, say $B_k:=A_{k0}$. 

Clearly, if $A_k$ satisfies a spectral condition as stated in the lemma, 
then $B_k$ does so. By \rf[wAk], it is also clear that the 
self-adjoint operators $B_0,\ldots,B_m$ strongly commute. 
Let $e(\lambda_0,\ldots,\lambda_m)$ denote the joint spectral 
projections. Take 
$\D_0:=\cup_{l\in\N}\big(e(l,\ldots,l)-e(-l,\ldots,-l)\big)\Hh$. 
Then $\D_0$ is an invariant core for each $B_j$
and $\D:=\Lin\{\eta_n\,;\,n\in\Z,\,\eta\in\D_0\}$ 
is an invariant core for each $A_j$. This completes the proof. 
\hfill$\Box$
\begin{thl}                           \label{L3}
Let $\epsilon\in\{\pm 1\}$. 
Assume that  $z$ is a closed operator on a Hilbert space $\Hh$.
Then we have $\D(zz^\ast)=\D(z^\ast z)$, this domain is dense in $\Hh$ 
and the relation 
\begin{equation}                                         \label{qhyp}
z^\ast z -q^2 z z^\ast = \epsilon(1-q^2)
\end{equation} 
holds if and only if $z$ is unitarily equivalent to an orthogonal direct 
sum of operators of the following form.
\begin{description}
\item[$\epsilon=1$:]
  \begin{description}   
   \item[$(I)$] 
     $z\eta_n = (1-q^{2(n+1)})^{1/2} \eta_{n+1}$ on 
     $\Hh={\oplus}^{\infty}_{n=0}\Hh_n$, 
     where each $\Hh_n$ is the same Hilbert space $\Hh_0$.
    \item[$(II)_A$] 
     $z\eta_n =(1+q^{2(n+1)} A^2)^{1/2}\eta_{n+1}$
     on $\Hh={\oplus}^{\infty}_{n=-\infty}\Hh_n$, $\Hh_n=\Hh_0$,
     where $A$ denotes a self-adjoint operator on 
     a Hilbert space 
     $\Hh_0$ such that $\sigma(A) \subseteq [q,1]$ 
and either $q$ or $1$ is not an eigenvalue of $A$.
   \item[$(III)_u$] $z=u$, where $u$ is a unitary
     operator on $\Hh$.
  \end{description}
\item[$\epsilon=-1$:]
  \begin{description}   
    \item[$(IV)$] $z\eta_n = (q^{-2n}-1)^{1/2} \eta_{n-1}$ 
    on $\Hh={\oplus}^{\infty}_{n=0}\Hh_n$, 
    where each $\Hh_n$ is $\Hh_0$.
  \end{description}
\end{description}
\end{thl}
The proof of Lemma \ref{L3} can be found in \cite{KW}. 
A version of this lemma appears also in \cite{CGP}, where  irreducible 
$\ast$-rep\-re\-sen\-ta\-tions of \rf[qhyp] are discussed.

\section{Cross product algebras related to 
                      the quantum motion group 
                       {\mathversion{bold}$\mathrm{E}_q(2)$} }
                                                     \label{sec-Eq2}
%
\subsection{Definitions}
                                                       \label{sec-def}
The coordinate Hopf $\ast$-al\-ge\-bra 
of the quantum motion group $\mathrm{E}_q(2)$ 
is the complex unital 
$\ast$-al\-ge\-bra $\cO(\mathrm{E}_q(2))$ with generators 
$v$, $v^\ast $, $n$, $n^\ast $ and defining relations 
\begin{equation}                                   \label{rnrel}
v^\ast v=vv^\ast=1,\quad nn^\ast=n^\ast n,\quad nv=q vn.
\end{equation}
The Hopf algebra structure is given by 
\begin{align}       \label{com}
&\Delta(v)=v\otimes v,\quad \Delta(n)=v\otimes n+n\otimes v^\ast,\\
&\varepsilon (v)=1, \quad \varepsilon(n)=0,\quad S(v)= v^\ast,\quad 
S(n)= -qn.                                         \label{com1}
\end{align}
Let $\cO(\C_q)$ denote the complex unital 
$\ast$-al\-ge\-bra with a single generator $z$ 
satisfying
\begin{equation}                                        \label{zrel}
z^\ast z=q^2zz^\ast.
\end{equation}
We call $\cO(\C_q)$ the coordinate $\ast$-al\-ge\-bra of the quantum complex 
plane. The $\ast$-al\-ge\-bra $\cO(\C_q)$ is a left and a right  
$\cO(\mathrm{E}_q(2))$-comodule 
$\ast$-al\-ge\-bra with left coaction $\varphi_L$  and right 
coaction $\varphi_R$  determined by
\begin{equation}                                      \label{zcoact}
\varphi_L(z)=v^2\otimes z+vn\otimes 1,\quad
\varphi_R(z)= 1\otimes vn^\ast  + z\otimes v^2.
\end{equation}

Using (\ref{rnrel}), we see that the map ${z{\rightarrow} vn}$ extends  
to a $\ast$-iso\-mor\-phism of $\cO(\C_q)$ onto the $\ast$-sub\-al\-ge\-bra of 
$\cO(\mathrm{E}_q(2))$ generated by $vn$. 
By (\ref{com}), \rf[com1] and (\ref{zcoact}), 
this $\ast$-iso\-mor\-phism intertwines the left coaction of 
$\cO(\mathrm{E}_q(2))$ and 
the comultiplication. Thus we can consider $\cO(\C_q)$ as a 
left $\cO(\mathrm{E}_q(2))$-comodule $\ast$-sub\-al\-ge\-bra of 
$\cO(\mathrm{E}_q(2))$ by 
identifying $z$ with $vn$. Analogously, we can consider $\cO(\C_q)$ as a 
right $\cO(\mathrm{E}_q(2))$-comodule $\ast$-sub\-al\-ge\-bra of 
$\cO(\mathrm{E}_q(2))$ by 
identifying $z$ with $vn^\ast $.

The quantized universal enveloping algebra $\U_q(\mathrm{e}_2)$ of 
the quantum motion group  is generated by $E$, $F$, $K$ and $K^{-1}$ 
with relations 
$$
KK^{-1} = K^{-1}K = 1,\quad KF = q FK,\quad KE = q^{-1}EK,
\quad EF = FE.
$$
It is a Hopf $\ast$-al\-ge\-bra with involution and Hopf algebra 
structure given by
\begin{align*}
&\mbox{ }\hspace{140pt} K^\ast=K,\quad E^\ast=F,\hspace{140pt} \mbox{ }\\
&\Delta(K)=K\otimes K, \ \ \Delta(E)=E\otimes K+K^{-1}\otimes E,\ \ 
\Delta(F)=F\otimes K+K^{-1}\otimes F,\\
&\varepsilon(K)=1, \ \ \varepsilon(E)=\varepsilon(F)=0, \ \ 
 S(K)=K^{-1},\ \   S(E)=-q^{-1}E, \ \  S(F)=-qF.
\end{align*}
There is a dual pairing of Hopf $\ast$-al\-ge\-bras
$\langle \cdot,\!\cdot\rangle : 
\U_q(\mathrm{e}_2)\times \cO(\mathrm{E}_q(2))\rightarrow \C$
which is zero on all pairs of generators except 
$$
\langle E,n^\ast\rangle =-q^{-1},\quad \langle F,n\rangle =1,\quad
\langle K,v\rangle =q^{1/2},\quad \langle K,v^\ast\rangle =q^{-1/2}.
$$

With the coaction induced by the comultiplication, 
$\Eq$ becomes a left and right $\Eq$-comodule $\ast$-al\-ge\-bra. 
By \rf[rechts] and \rf[ang], $\Eq$ is a right and a left 
$\ue$-mod\-ule $\ast$-al\-ge\-bra. Simple computations show that the 
right action $\anf$ and the left action $\ang$ are given 
by 
\begin{align*}
 &n^\ast\anf E=-q^{-1} v,& &n\anf E=v\anf E=v^\ast\anf E=0, \\
 &n\anf F=v^\ast,& &n^\ast\anf F=v\anf F=v^\ast\anf F=0,\\
 &n\anf K=q^{1/2} n,& &n^\ast\anf K=q^{-1/2} n^\ast,\quad\ 
   v\anf K=q^{1/2} v,\quad\ v^\ast\anf K=q ^{-1/2} v^\ast;\\[-24pt]
\end{align*}
\begin{align*}
 &E\ang n^\ast = -q^{-1} v^\ast,& &E\ang n=E\ang v=E\ang v^\ast=0,\\
 &F\ang n=v,&  &F\ang n^\ast=F\ang v=F\ang v^\ast=0,\\
 &K\ang n=q^{-1/2}n,&  &K\ang n^\ast=q^{1/2}n^\ast,\quad\
K\ang v=q^{1/2}v,\quad\ K\ang v^\ast=q^{-1/2}v^\ast.
\end{align*}
Similarly, by \rf[rechts], \rf[ang] and \rf[zcoact], 
$\oc$ is a left and a right $\ue$-mod\-ule $\ast$-al\-ge\-bra 
with left and right action determined by
$$
z\anf K=q z,\ z^\ast\anf K=q^{-1} z^\ast,\ 
z\anf E=0,\ z^\ast\anf E=-q^{-3/2},\  z\anf F=q^{-1/2},\ z^\ast\anf F=0;
$$  
$$
K\ang z=q z,\ K\ang z^\ast=q^{-1}z^\ast,\ 
E\ang z=-q^{-3/2},\ E\ang z^\ast=0,\ 
F\ang z=0,\ F\ang z^\ast= q^{-1/2}. 
$$

From \rf[cross1], \rf[oma] and above formulas, 
we derive the following cross commutation relations in  
the corresponding cross product algebras. 
\begin{align*}
\U_q(\mathrm{e}&_2)\lti\cO(\mathrm{E}_q(2)):  & 
 vF&=q^{1/2} Fv, &  vE&=q^{1/2} Ev, &   vK&=q^{1/2} Kv,\\
& & v^\ast F&=q^{-1/2} Fv^\ast,  &    v^\ast E&=q^{-1/2} Ev^\ast,
 &     v^\ast K&=q^{-1/2}Kv^\ast,\\
 & 
\makebox[0cm][l]{$
  nF=q^{1/2}Fn+K^{-1}v^\ast,\quad nE=q^{1/2}En, \quad  nK=q^{1/2} Kn,$}
& & & & & &\\
 & 
\makebox[0cm][l]{$
 n^\ast F=q^{-1/2} Fn^\ast,\quad 
n^\ast E=q^{-1/2} En^\ast-q^{-1}K^{-1}v,\quad n^\ast K=q^{-1/2}Kn^\ast .$}
& & & & & &\\[-24pt]
\end{align*}
\begin{align*}
\cO(\mathrm{E}&_q(2))\rti\ue: &
Fv&=q^{-1/2}vF, &  Ev&=q^{-1/2}vE, &  Kv&=q^{1/2}vK, \\
& & Fv^\ast&=q^{1/2}v^\ast F, &  
Ev^\ast&=q^{1/2}v^\ast E, &  Kv^\ast &=q^{-1/2}v^\ast K,\\
&\makebox[0cm][l]{$Fn=q^{1/2}nF+vK,\quad 
En=q^{1/2}nE,\quad  Kn=q^{-1/2}nK,$} \\
& \makebox[0cm][l]{$Fn^\ast=q^{-1/2}n^\ast F,\quad
 En^\ast=q^{-1/2}n^\ast E-q^{-1}v^\ast K,
\quad Kn^\ast=q^{1/2}n^\ast K. $} \\[-24pt]
\end{align*}
\begin{align*}
\U_q(\mathrm{e}_2)\lti\cO(\C_q): \qquad
zK=qKz, \quad zE=qEz, \quad zF=qFz + q^{-1/2} K^{-1}, \\
\qquad 
z^\ast K=q^{-1}Kz^\ast, \quad z^\ast E=q^{-1}Ez^\ast - q^{-3/2}K^{-1},
\quad z^\ast F=q^{-1}Fz^\ast.\\[-24pt]
\end{align*}
\begin{align*}
\cO(\C_q)\rti\U_q(\mathrm{e}_2): \qquad
Kz=qzK,\quad Ez=q^{-1}zE-q^{3/2}K,\quad Fz=q^{-1}zF,\\
\qquad
Kz^\ast=q^{-1}z^\ast K,\quad Ez^\ast=qz^\ast E,\quad 
Fz^\ast=qz^\ast F+q^{-1/2}K.
\end{align*}

The next lemma shows that we can restrict ourselves to the right 
versions of these cross product algebras. 
It is proved by direct computations. 
\begin{thl}
There is a $\ast$-iso\-mor\-phism 
$$
\theta\,:\,\U_q(\mathrm{e}_2)\lti\cO(\mathrm{E}_q(2))\longrightarrow
\cO(\mathrm{E}_q(2))\rti\ue
$$ 
determined by $\theta (v)=v$, $\theta (n)=n^\ast$,   
$\theta(K)=K^{-1}$ and $\theta (E)=F$.  

There is a $\ast$-iso\-mor\-phism 
$$
\psi\,:\,\U_q(\mathrm{e}_2)\lti\cO(\C_q)\longrightarrow
\cO(\C_q)\rti\U_q(\mathrm{e}_2)
$$ 
determined by $\psi (z)=z$, $\psi(K)=K^{-1}$ and $\psi (E)=F$. 

The inverse isomorphisms $\theta^{-1}$ 
and $\psi^{-1}$ are given by the same formulas.
\end{thl}

Let $\U_0$ denote the subalgebra of $\U_q(\mathrm{e}_2)$ generated by 
the unit element and the linear span $\T_0$ of the elements 
$$
X:=q^{1/2} FK,\qquad Y:=-q^{3/2}EK.
$$
Clearly, 
$YX=q^2XY$ and $X^\ast =-q^{-2}Y$.
In particular, $\U_0$ 	is a $\ast$-al\-ge\-bra. Since $\C\cdot 1+\T_0$ is 
a right coideal of $\U_q(\mathrm{e}_2),\T_0$ is the quantum tangent space of 
a left-covariant first order differential calculus on 
$\cO(\mathrm{E}_q(2))$ \cite[Proposition 14.5]{KS}. 
It can be shown that this calculus induces a differential 
$\ast$-cal\-cu\-lus on the 
$\ast$-sub\-al\-ge\-bra $\cO(\C_q)$ such that $X$ and $Y$ 
can be considered as partial derivatives. 
We shall not carry out the details because we are  interested 
in the $\ast$-al\-ge\-bra $\U_0\lti\cO(\C_q)$ only. 

The $\ast$-al\-ge\-bra $\U_0\lti\cO(\C_q)$ is the $\ast$-sub\-al\-ge\-bra of 
$\ue\lti\cO(\C_q)$ generated by $\U_0$ and $\cO(\C_q)$ or, equivalently,  
the $\ast$-al\-ge\-bra with generators $X$, $X^\ast$, $z$, $z^\ast$ 
and defining relations  
\begin{align}  
\U_0\lti\cO(\C_q): \qquad                                                 
z^\ast z&=q^2zz^\ast, &  X^\ast X&=q^2XX^\ast,           \label{qp1}\\
zX&=q^2Xz+1, &  zX^\ast&=q^2X^\ast z,                    \label{qp2}\\
z^\ast X&=q^{-2}Xz^\ast, &                                   \label{qp3}
z^\ast X^\ast&=q^{-2}X^\ast z^\ast-q^{-2}.\qquad
\end{align}
%
%
\subsection{Representations of the {\mathversion{bold}$\ast$}-al\-ge\-bra 
             {\mathversion{bold}$ \U_0\!\ltimes\!\cO(\C_q)$}}
                                                         \label{sec-U0Cq}
We set $\gamma =(1-q^2)^{-1}$ and define $N=zX - \gamma$.  
In the $\ast$-al\-ge\-bra\ \, ${\U_0{\lti}\cO(\C_q)}$, the element $N$  
satisfies the following relations
\begin{equation}                                          \label{relnz}
zN=q^2Nz,\quad z^\ast N=Nz^\ast,\quad N^\ast N=NN^\ast, \quad
z^\ast zN^\ast N=N^\ast Nz^\ast z.
\end{equation}
These four equations follow immediately from 
(\ref{qp1})--(\ref{qp3}). 
As a sample, we verify the relation $N^\ast N=NN^\ast$. 
Indeed, from
$$
X^\ast z^\ast z X=q^2 X^\ast z z^\ast X=X^\ast z Xz^\ast
=q^{-2} z X^\ast Xz^\ast=zXX^\ast z^\ast, 
$$
we conclude 
\begin{align*}
N^\ast N\hspace{-1.23pt}=\hspace{-1.23pt} X^\ast z^\ast zX
\hspace{-1.23pt}-\hspace{-1.23pt}\gamma z X\hspace{-1.23pt}-\hspace{-1.23pt}
\gamma X^\ast z^\ast \hspace{-1.23pt}+\hspace{-1.23pt}\gamma^2
\hspace{-1.23pt}=\hspace{-1.23pt}zXX^\ast z^\ast
\hspace{-1.23pt}-\hspace{-1.23pt}\gamma zX
\hspace{-1.23pt}-\hspace{-1.23pt}\gamma X^\ast z^\ast 
\hspace{-1.23pt}+\hspace{-1.23pt}\gamma^2
\hspace{-1.23pt}=\hspace{-1.23pt}NN^\ast.
\end{align*}

Now suppose we are given  a $\ast$-rep\-re\-sen\-ta\-tion 
of the $\ast$-al\-ge\-bra $\U_0\lti\cO(\C_q)$ on a 
Hilbert space $\Hh$. 
As explained in the introduction, we assume 
that $N$ is a normal operator and that 
the self-adjoint operators $z^\ast z$ and $N^\ast N$
strongly commute. 

We claim that $\ker z = \ker z^\ast = \{0\}$.
To see this, observe that 
$z^\ast z=q^2 zz^\ast$ yields  
$\ker z=\ker z^\ast$. Let $\eta\in\ker z$. 
Then, by (\ref{qp2}), 
$$
\parallel \eta\parallel^2=\langle (zX-q^2 Xz)\eta,\eta\rangle 
=\langle zX\eta,\eta\rangle =\langle X \eta,z^\ast\eta\rangle =0
$$
so that $\eta=0$. Thus, $\ker z=\{0\}$. 

Let $z=w|z|$ be the polar decomposition of the closed operator $z$. 
As $\ker z=\ker z^\ast=\{0\}$, $w$ is unitary. 
The operator relation $z^\ast z=q^2 zz^\ast$ is equivalent to 
$|z|^2=q^2 w|z|^2 w^\ast$ and so to $|z|=q w|z|w^\ast$. 
Since $\ker |z|=\{0\}$, it follows from Lemma \ref{1} that there is 
a Hilbert space $\Hh_0$ and a self-adjoint operator $z_0$ on $\Hh_0$ 
satisfying $\sigma(z_0)\sqsubseteq (q,1]$ such that
$$
\Hh=\displaystyle\mathop{\oplus}^\infty_{n=-\infty}\Hh_n,\ 
\Hh_n=\Hh_0, \quad   w\eta_n=\eta_{n-1},\quad 
|z|\eta_n=q^{-n} z_0\eta_n.
$$

As $z^\ast z$ and $N^\ast N$ strongly commute,  
$|z|$ and $|N|$ do so. Consequently, $|N|$ commutes 
with the spectral projections of $|z|$. 
Since $\sigma (z_0)\sqsubseteq (q,1]$, it follows that $|N|$ 
leaves each Hilbert space $\Hh_n$ invariant, 
so there are self-adjoint operators $N_n$
on $\Hh_0$ such that 
$$
|N|\eta_n=N_n \eta_n,\quad  \eta_n\in\Hh_n.
$$

An application of \rf[relnz] yields 
$zN^\ast N=q^2 N^\ast Nz$ which entails the 
operator identity $w|z||N|^2=q^2|N|^2 w|z|$. 
Using above commutation relations, 
we get $|z|w|N|^2=q^2|z||N|^2 w$. 
Since $\ker |z| = \{0\}$, 
it follows that  $w |N|^2=q^2 |N|^2w$. Hence $w|N|=q|N| w$ 
because $w$ is unitary. This in turn gives
$N_n\eta_{n-1}=w|N|w^\ast \eta_n= q N_{n-1}\eta_{n-1}$ and so 
$N_n=q^n N_0$. Accordingly, $|N|\eta_n=q^nN_0\eta_n$.

From \rf[relnz], $zN^\ast N=q^2 N^\ast N z$ 
and  $N^\ast N z^\ast =q^2 z^\ast N^\ast N$. 
We assume that these identities hold in strong sense. 
Since $\ker N= \ker N^\ast N$, we see that 
$\ker N$ is reducing. 
Hence the representation decomposes into a direct sum of 
two representations corresponding to the cases $N=0$ and $\ker N=\{0\}$. 
We treat the two cases separately.

\mn
{\it Case I.} $N=0$. 

\nopagebreak

This means that $X=(1 - q^2)^{-1} z^{-1}$. 
Using the fact that $z^\ast z=q^2zz^\ast$, one immediately 
verifies that all relations (\ref{qp1})--(\ref{qp3}) are fulfilled.

\mn
{\it Case II.} $\ker N=\{0\}$. 

\nopagebreak

Using this additional condition, we continue the above operator-theoretic 
manipulations. Let $N = u|N|$ be the polar decomposition of $N$. 
Since $\ker N = \{0\}$ and $N$ is normal, $u$ is unitary 
and $u |N|=|N|u$. 
Applying (\ref{relnz}) 
and the assumption that $|z|$ and $|N|$ strongly commute, we obtain
$$
|z|^2 u |N|=z^\ast z N=q^2 Nz^\ast z=q^2 u |N| |z|^2=q^2 u |z|^2|N|.
$$
As $\ker |N|=\{0\}$, above equation  implies $|z|^2 u=q^2u |z|^2$. 
The operator $u$ is unitary, thus $|z| u=qu|z|$. 
Since $|z|w=qw|z|$, $w^\ast u$ commutes with the self-adjoint 
operator $|z|$ and hence with the spectral projections of $|z|$. 
Therefore, $w^\ast u$ leaves each space $\Hh_n$ invariant. 
Hence there are unitary operators $u_n$ on $\Hh_0$, $n\in\Z$, such 
that $w^\ast u\eta_n=u_n\eta_n$ for $\eta_n \in\Hh_n$. Accordingly, 
$u\eta_n=u_n\eta_{n-1}$. 
Using the relations  $|N| w^\ast=q w^\ast |N|$, $u |N|=|N|u$, 
$|z| u=qu|z|$ and the fact that $|z|$ and $|N|$ 
strongly commute, we derive 
%
$$
|z||N|w^\ast u=q|z|w^\ast u|N|=qz^\ast N=qNz^\ast= qu |N| |z| w^\ast
=|z| |N|uw^\ast .
$$
Thus, $w^\ast u= uw^\ast$ since $\ker |z|=\ker |N|=\{0\}$. This gives 
$u_n\eta_{n}=uw^\ast\eta_{n-1}=w^\ast u\eta_{n-1}= u_{n-1}\eta_n$, 
so $u_n=u_{n-1}$. Hence $u_n=u_0$  for all $n\in\Z$. Employing
the relations $w^\ast u|z|=|z|w^\ast u$ and $u |N|=|N|u$,  
we derive  $u_0z_0=z_0u_0$ and $u_0N_0u^\ast_0=q^{-1} N_0$. 
Since $|N|$ and $|z|$ strongly commute, $N_0z_0=z_0N_0$. 
As $|N|$ is a positive self-adjoint operator with trivial kernel, 
so is $N_0$. 
Inserting $z^{-1}=|z|^{-1}w^\ast$ and $N=u|N|$ into $X=z^{-1}(N+\gamma)$, 
we can express the operator $X$ (and its adjoint $X^\ast$) 
in terms of $u_0$, $z_0$ and $N_0$. 

Summarizing, we have in Case II 
\begin{align}                                               
X\eta_n&=q^{2n} z^{-1}_0 u_0 N_0 \eta_n 
+ (1 - q^2)^{-1} q^{n+1} z^{-1}_0 \eta_{n+1},           \label{opx}\\
X^\ast \eta_n &= q^{2n} z^{-1}_0 N_0 u^\ast_0 \eta_n
+ (1 - q^2)^{-1} q^n z^{-1}_0 \eta_{n-1},               \label{opxstav}\\
z\eta_n&=q^{-n} z_0 \eta_{n-1},\quad 
z^\ast \eta_n = q^{1-n} z_0 \eta_{n+1}                  \label{opz}
\end{align}
on the Hilbert space 
$\Hh=\oplus^\infty_{n=-\infty} \Hh_n$, $\Hh_n=\Hh_0$, 
where $z_0$, $N_0$ are self-adjoint operators and $u_0$ is 
a unitary operator on $\Hh_0$ such that 
$\sigma (z_0)\sqsubseteq (q,1]$, $N_0>0$, and
\begin{equation}                                        \label{conznv}
z_0N_0=N_0z_0,\quad u_0z_0u_0^\ast=z_0,\quad u_0N_0u^\ast_0=q^{-1} N_0
\end{equation}
Conversely, if the latter is satisfied, 
then formulas (\ref{opx})--(\ref{opz}) define a
$\ast$-rep\-re\-sen\-ta\-tion of $\U_0\lti\cO (\C_q)$.
\mn

The $\ast$-rep\-re\-sen\-ta\-tions of the relations \rf[conznv] are 
described by Lemma \ref{1}. 
Inserting the expressions for $z_0$, $N_0$ and $u_0$ 
from Lemma \ref{1} into Equations \rf[opx]--\rf[opz]   
and renaming suitable, we obtain 
the following list of $\ast$-rep\-re\-sen\-ta\-tion of 
$\U_0\lti\cO(\C_q)$. 
\begin{align*}
(I)_A:&  &z \eta_n&=q^{-n} A\eta_{n-1},\quad  
         z^\ast \eta_n =q^{-(n+1)}A\eta_{n+1},\\
& & X\eta_n&= (1 - q^2)^{-1} q^{(n+1)} A^{-1} \eta_{n+1},\\
& & X^\ast\eta_n&=(1 - q^2)^{-1}q^{n} A^{-1} \eta_{n-1}
\quad \mbox{on} \ \,
\Hh = {\displaystyle\mathop{\oplus}^{\infty}_{n=-\infty}} \Hh_n,\ 
\Hh_n=\K.\\[-24pt]
\end{align*}
\begin{align*}
(II)_{A,B}:& &z\eta_{nk}&=q^{-n} A\eta_{n-1,k},\quad  
      z^\ast \eta_{nk} =q^{-(n+1)}A \eta_{n+1,k},\\
  & & X\eta_{nk} &= q^{2n-k} A^{-1} B\eta_{n,k-1} +(1 - q^2)^{-1}      
        q^{n+1} A^{-1} \eta_{n+1,k},\\
  & & X^\ast \eta_{nk} &= q^{2n-k-1} A^{-1} B\eta_{n,k+1} 
            +(1 - q^2)^{-1} q^n A^{-1} \eta_{n-1,k}\\
 & & &\hspace{4.5cm} \mbox{on} \ \,
\Hh = {\displaystyle\mathop{\oplus}^{\infty}_{n,k=-\infty}} \Hh_{nk},\ 
\Hh_{nk}=\K.
\end{align*}
The parameters $A$ and $B$ denote self-adjoint operators on the 
Hilbert space $\K$ such that 
$\sigma (A)\sqsubseteq (q,1]$, $\sigma (B)\sqsubseteq (q,1]$,  
and, in the case $(II)_{A,B}$,  $AB=BA$. 
Rep\-re\-sen\-ta\-tions 
labeled by different sets of parameters (within unitary equivalence)
or belonging to 
different series are not unitarily equivalent. 
A representation is irreducible 
if and only if $\K$ is isomorphic to the one-dimensional 
Hilbert space $\C$. In this case, we can regard $A$ and $B$ as 
real numbers of the interval $(q,1]$.
%
%
\subsection{Representations of the {\mathversion{bold}$\ast$}-al\-ge\-bra 
            {\mathversion{bold}$\U_q(\mathrm{e}_2)\!\ltimes\!\cO(\C_q)$}}
                                                  \label{rep-Cq}
Suppose we have a $\ast$-rep\-re\-sen\-ta\-tion of the 
$\ast$-al\-ge\-bra 
$\U_q(\mathrm{e}_2){\lti}\cO(\C_q)$ on a Hilbert space $\Hh$. 
Then the considerations of 
the preceding subsection apply to the restriction of the representation 
to the subalgebra $\U_0\lti\cO(\C_q)$. We freely use the facts and notations 
set up therein.

The new ingredient in the larger algebra 
$\U_q(\mathrm{e}_2)\lti\cO(\C_q)$ is 
the invertible generator $K$ satisfying the relations
\begin{equation}                                         \label{K}
zK=qKz,\quad z^\ast K=q^{-1}Kz^\ast,\quad XK=q^{-1} KX,
\quad X^\ast K=qKX^\ast .
\end{equation}
Again, we assume that these relations hold in 
strong operator-theoretic sense. 

\mn
{\it Case I.} 
Let $z$, $X$ and $X^\ast$ be given as described in $(I)_A$.
By \rf[K], $K$ commutes with $z^\ast z$ and 
hence with the spectral projections of $z^\ast z$. 
Since $z^\ast z\eta_n=q^{-2n} A^2\eta_n$ and 
$\sigma(A)\sqsubseteq (q,1]$, we conclude that $K$ leaves each space 
$\Hh_n$ invariant, that is, there are operators $K_n$ on $\Hh_n$, $n\in\Z$, 
such that $K\eta_n=K_n\eta_n$. Inserting the expressions of $z$ and $K$ 
into the relations  $zz^\ast K=Kz^\ast z$ and $zK=qKz$
gives  $AK_n=K_n A$ and $K_n=q K_{n-1}$. Setting $H:=K_0$, we can write 
$K_n=q^n H$, where $H$ is an invertible self-adjoint operator on $\K$ 
commuting with $A$.  
Finally, $N\equiv zX - \gamma=0$ and so 
$F=q^{-1/2} XK^{-1}= q^{-1/2} \gamma z^{-1}K^{-1}$.
This determines the actions of the operators $z$, $z^\ast$, $K$, $F$ 
and $E=F^\ast$ completely.

\mn
{\it Case II.}
Suppose that the representation of the  operators $z$, $X$ and $X^\ast$ 
takes the form described in $(II)_{A,B}$. 
As in the preceding paragraph, it follows from \rf[K] that $K$ is given 
on $\Hh_n=\oplus^{\infty}_{k=-\infty}\Hh_{nk}$ 
by $K\zeta_n=q^n K_0\zeta_n$, $\zeta_n\in\Hh_n$, where $K_0$ is 
an invertible self-adjoint operator acting on $\Hh_0$.  
Relations \rf[K] and the definition of $N$ yield $NK=KN$ and 
$N^\ast K=K N^\ast$. Thus we can assume that the self-adjoint 
operators $K$ and $|N|$ strongly commute. Observe that $|N|$ acts on 
$\Hh_0$ by $|N|\eta_{0k}=q^{-k}B\eta_{0k}$ and that  
$\sigma (B)\sqsubseteq (q,1]$. Since $K$ commutes with the spectral 
projections of $|N|$, it leaves the Hilbert spaces $\Hh_{0k}$ invariant. 
Hence there exist invertible self-adjoint operators 
$K_{0k}$ on $\Hh_{0k}$ such that 
$K_0\eta_{0k}=K_{0k}\eta_{0k}$ and $B K_{0k}=K_{0k} B$. 
From the last subsection, $N\eta_{nk}=q^{n-k}B\eta_{n-1,k-1}$. 
Applying $NK=KN$ to vectors $\eta_{nk}\in\Hh_{nk}$ gives 
$q^{2n-k}K_{0k}B\eta_{n-1,k-1}=q^{2n-k-1}K_{0,k-1}B\eta_{n-1,k-1}$, 
hence $K_{0k}=q^{-1}K_{0,k-1}$. 
Denoting $K_{00}$ by $H$, we get $K\eta_{nk}=q^{n-k}H\eta_{nk}$, 
where $H$ is an invertible self-adjoint operator on $\Hh_{00}$
commuting with $B$. Moreover, it follows from 
$z^\ast zK=Kz^\ast z$ that $H$ commutes with $A$.  
The actions of $F$ and $E$ on $\Hh$ are obtained  by computing 
$F=q^{-1/2}z^{-1}(N+\gamma)K^{-1}$ and $F^\ast=E$. 

Summarizing, we have obtained the following series of 
$\ast$-rep\-re\-sen\-ta\-tions of $\U_q(\mathrm{e}_2)\lti\cO(\C_q)$.
\begin{align*}
(I)_{A,H}:&  &z \eta_n&=q^{-n} A\eta_{n-1},\quad  
         z^\ast \eta_n =q^{-(n+1)}A\eta_{n+1},\\
& & F\eta_n&=q^{1/2} (1 - q^2)^{-1} A^{-1} H^{-1}\eta_{n+1},\\
& & E\eta_n&=q^{1/2} (1 - q^2)^{-1} A^{-1} H^{-1}\eta_{n-1},\\
& & K\eta_n&=q^{n}H\eta_{n}
\quad \mbox{on} \ \,
\Hh = {\displaystyle\mathop{\oplus}^{\infty}_{n=-\infty}} \Hh_n,\ 
\Hh_n=\K.\hspace{2.5cm}\\[-24pt]
\end{align*}
\begin{align*}
(II)_{A,B,H}:\qquad &z\eta_{nk}=q^{-n} A\eta_{n-1,k},\quad  
      z^\ast \eta_{nk} =q^{-(n+1)}A \eta_{n+1,k},\\
   F\eta_{nk} &= q^{-1/2}q^{n} A^{-1} BH^{-1}\eta_{n,k-1} 
   +q^{1/2}(1 - q^2)^{-1}q^{k} A^{-1}H^{-1} \eta_{n+1,k},\\
   E\eta_{nk} &= q^{-1/2}q^{n} A^{-1} BH^{-1}\eta_{n,k+1} 
   +q^{1/2}(1 - q^2)^{-1}q^{k} A^{-1}H^{-1} \eta_{n-1,k},\\
   K\eta_{nk}&=q^{n-k}H\eta_{nk}
\quad \mbox{on} \ \,
\Hh = {\displaystyle\mathop{\oplus}^{\infty}_{n,k=-\infty}} \Hh_{nk},\ 
\Hh_{nk}=\K.
\end{align*}
Here $A$, $B$ and $H$ are commuting self-adjoint operators on the 
Hilbert space $\K$ such that 
$\sigma (A)\sqsubseteq (q,1]$, $\sigma (B)\sqsubseteq (q,1]$,  and 
$H$ is invertible. 
Rep\-re\-sen\-ta\-tions 
labeled by different sets of parameters (within unitary equivalence)
or belonging to 
different series are not unitarily equivalent. 
A representation is irreducible if and only if $\K=\C$. 
Then $A$, $B$ and $H$ are real numbers such that 
$A,B\in(q,1]$ and $H\neq 0$. 
%
%
\subsection{Representations of the {\mathversion{bold}$\ast$}-al\-ge\-bra 
   {\mathversion{bold}$ \U_q(\mathrm{e}_2)\!\ltimes\!\cO(\mathrm{E}_q(2))$}}
                                                         \label{rep-Eq2}
Recall from Subsection \ref{sec-def} that 
$\oc$ becomes a 
$\ast$-sub\-al\-ge\-bra of                                        
$\Eq$ by identifying $z$ with $vn$. 
The relations from Subsection \ref{sec-def} show that   
the cross product algebra  
$\U_q(\mathrm{e}_2)\lti\cO(\mathrm{E}_q(2))$                             
can be described as the $\ast$-al\-ge\-bra 
generated by $z, E, F, K \in \U_q(\mathrm{e}_2)\lti\cO(\C_q)$ 
and an additional generator $v$ satisfying 
\begin{equation}                                          \label{addrel}
vv^\ast=v^\ast v=1,\  zv=q vz,\  vK=q^{1/2} Kv,\  
vE=q^{1/2} Ev,\  vF=q^{1/2} Fv.
\end{equation}
The element $n\in\cO(\mathrm{E}_q(2))$ is recovered by setting 
$n=v^\ast z$. 

For a $\ast$-rep\-re\-sen\-ta\-tions of 
$\U_q(\mathrm{e}_2)\lti\cO(\mathrm{E}_q(2))$, 
we apply the results of the preceding subsection to its restriction to the 
$\ast$-sub\-al\-ge\-bra $\U_q(\mathrm{e}_2)\lti\cO(\C_q)$. 
It only remains to determine the 
action of the additional operator $v$. 

\mn
{\it Case I.} 
Assume that the operators $z$, $E$, $F$, $K$ are given by 
the formulas of the series $(I)_{A,H}$. 
Let $z=w|z|$ be the polar decomposition of the closed operator $z$. 
Then $w$ and $|z|$ act on 
$\Hh = \oplus^{\infty}_{n=-\infty} \Hh_n$, 
$\Hh_n=\Hh_0$, by $w\eta_n=\eta_{n-1}$ and $|z|\eta_n=q^{-n}A\eta_n$, 
where $\eta_n\in\Hh_n$.
The relations $vv^\ast=v^\ast v=1$ and $zv=q vz$ yield
$z^\ast v=q vz^\ast$ and so $z^\ast zv=q^2vz^\ast z$. 
This implies $|z|v=qv|z|$ since $v$ is unitary. 
Note that $|z|w=qw|z|$. It follows that the unitary operator 
$w^\ast v$ commutes with $|z|$. As a consequence,  
$w^\ast v$ leaves each space $\Hh_n$ invariant. 
Hence there are unitary operators  $v_n$ on $\Hh_n$ 
such that $w^\ast v\eta_n=v_n\eta_n$.
Accordingly, $v\eta_n=v_{n} \eta_{n-1}$. 
Evaluating $|z|v=qv|z|$ and $zv=qvz$ on vectors $\eta_n\in\Hh_n$ 
gives $Av_{n}=v_{n}A$ and $Av_{n}=v_{n-1}A$, respectively. 
Combining the first equation with the second shows that 
$v_{n}=v_{n-1}$ since $A$ is invertible. Hence we can write 
$v\eta_n=v_0\eta_{n-1}$, where $v_0$ is a unitary operator on $\Hh_0$ 
commuting with $A$. 
Finally, by applying $vK=q^{1/2} Kv$ to vectors $\eta_n$, 
one gets $v_0H=q^{-1/2} Hv_0$ on $\Hh_0$. 
It can easily be checked that, if the preceding conditions on the operator 
$v$ are satisfied, then the identities \rf[addrel] hold. 

In conclusion, we have to determine operators $v_0$, $A$, $H$ 
on $\Hh_0$ satisfying 
\begin{equation}                                        \label{vAH}
Av_{0}=v_{0}A,\quad Hv_0=q^{1/2}v_0H ,\quad AH=HA,
\end{equation}
where $v_0$ is unitary, $H$ is an invertible self-adjoint operator, 
and $A$ is a self-adjoint operator such that 
$\sigma (A)\sqsubseteq (q,1]$.
Note that the subspaces of $\Hh_0$ where $H>0$ and $H<0$ 
are reducing. Considering the cases  $H>0$ and $H<0$ separately, 
we can apply Lemma \ref{1} to establish the 
$\ast$-rep\-re\-sen\-ta\-tions of the relations \rf[vAH].

\mn 
{\it Case II.} 
Suppose we are given operators $z$, $E$, $F$, $K$ as described by 
the formulas of the series $(II)_{A,B,H}$. 
As in Case I, one  shows by using 
$zv=q vz$ and $|z|v=qv|z|$  that $v$ maps 
$\Hh_n = \oplus^{\infty}_{j=-\infty} \Hh_{nj}$
into 
$\Hh_{n-1}=\oplus^{\infty}_{j=-\infty} \Hh_{n-1,j}$.
On the other hand, 
observe that $v$ commutes with $N=q^{1/2}zFK-(1-q^2)^{-1}$ and its 
adjoint. This yields $v|N|=|N|v$ since  $v$ is unitary. 
The action of $|N|$ on $\Hh_{nk}$
is given by $|N|\eta_{nk}=q^{n-k}B\eta_{nk}$, 
where $B$ is a self-adjoint operator on $\Hh_{00}$ such that 
$\sigma (B)\sqsubseteq (q,1]$. 
Since $v$ commutes with the spectral projections of $|N|$, 
we conclude that $v$ leaves each Hilbert space 
$\oplus^{\infty}_{j=-\infty} \Hh_{n+j,k+j}$ 
invariant. 
Combining these facts, it follows that $v$ maps $\Hh_{nk}$ into 
$\Hh_{n-1,k-1}$. Write $v\eta_{nk}=v_{nk}\eta_{n-1,k-1}$, where 
$v_{nk}$ denotes a unitary operator acting on $\Hh_{nk}=\Hh_{00}$. 
Applying $|z|v=q v|z|$ and $|N|v=v|N|$ to vectors $\eta_{nk}\in\Hh_{nk}$ 
gives $Av_{nk}=v_{nk}A$ and $Bv_{nk}=v_{nk}B$, respectively. 
Similarly, the relations $zv=q vz$ and $Nv=vN$ lead to  
$Av_{nk}=v_{n-1,k}A$ and $Bv_{nk}=v_{n-1,k-1}B$, respectively.
This yields 
$v_{nk}A=v_{n-1,k}A$ and $v_{nk}B=v_{n-1,k-1}B$.
Since $A$ and $B$ are invertible, we have 
$v_{nk}=v_{n-1,k}$ and $v_{nk}=v_{n-1,k-1}$. Hence 
all $v_{nk}$ are equal, say to $v_{00}$. 
Inserting the expression obtained for $v$ into 
$vK\eta_{nk}=q^{1/2} Kv\eta_{nk}$, we get 
$q^{n-k}v_{00}H\eta_{n-1,k-1}=q^{1/2}q^{n-k}H v_{00}\eta_{n-1,k-1}$ 
so that $v_{00}H=q^{1/2} Hv_{00}$.

Gathering the facts of Case II together shows that 
the representation is determined by 
operators $v_{00}$, $A$, $B$, $H$ 
on $\Hh_{00}$ satisfying 
\begin{align*}                                       
&Av_{00}=v_{00}A,\quad v_{00}H=q^{1/2} Hv_{00},\quad AH=HA,\\
& Bv_{00}=v_{00}B,\quad BH=BH ,\quad BA=AB, \notag
\end{align*}
where $v_{00}$ is unitary and 
$A$, $B$, $H$ are self-adjoint operators such that 
$\sigma (A)\sqsubseteq (q,1]$, $\sigma (B)\sqsubseteq (q,1]$, and
$H$ is invertible. 
Similarly to Case I, we employ Lemma \ref{1} to describe the 
representations of the above relations.

After a slight change of notation, 
the preceding discussion leads to the following
$\ast$-rep\-re\-sen\-ta\-tions of the cross product algebra 
$\U_q(\mathrm{e}_2)\lti\cO(\mathrm{E}_q(2))$.
\begin{align*}
(I)_{A,H,\epsilon}:&\mbox{ }  &v \eta_{mk}&=\eta_{m-1,k-1},\quad  
        n \eta_{mk} =q^{-m}A\eta_{m,k+1},\\
& & F\eta_{mk}&= 
(1 - q^2)^{-1}q^{(k+1)/2}\epsilon A^{-1} H^{-1}\eta_{m+1,k},\\
& & E\eta_{mk}&= 
(1 - q^2)^{-1}q^{(k+1)/2}\epsilon  A^{-1} H^{-1}\eta_{m-1,k},\\
& & K\eta_{mk}&=q^{m-k/2}\epsilon H\eta_{mk}
\quad \mbox{on} \ \,
\Hh = {\displaystyle\mathop{\oplus}^{\infty}_{m,k=-\infty}} \Hh_{mk},\ 
\Hh_{mk}=\K.\hspace{1.5cm}\\[-24pt]
\end{align*}
\begin{align*}
{(II)_{A,}}&_{B,H,\epsilon}: \qquad v\eta_{mkl}=\eta_{m-1,k-1,l+1},\quad  
      n \eta_{mkl} =q^{-m}A \eta_{m,k+1,l-1},\\
   F\eta_{mkl} &= q^{m+(l-1)/2}\epsilon A^{-1} BH^{-1}\eta_{m,k-1,l} 
   +(1 - q^2)^{-1}q^{k+(l+1)/2}\epsilon A^{-1}H^{-1} \eta_{m+1,kl},\\
   E\eta_{mkl} &= q^{m+(l-1)/2}\epsilon A^{-1} BH^{-1}\eta_{m,k+1,l} 
   +(1 - q^2)^{-1}q^{k+(l+1)/2}\epsilon A^{-1}H^{-1} \eta_{m-1,kl},\\
   K\eta_{mkl}&=q^{m-k-l/2}\epsilon H\eta_{mkl}
\quad \mbox{on} \ \,
\Hh = {\displaystyle\mathop{\oplus}^{\infty}_{m,k,l=-\infty}} \Hh_{mkl},\ 
\Hh_{mkl}=\K.
\end{align*}
The parameters $A$, $B$, $H$ denote  
commuting self-adjoint operators 
on a Hilbert space $\K$ such that 
$\sigma (A)\sqsubseteq (q,1]$, $\sigma (B)\sqsubseteq (q,1]$, 
$\sigma (H)\sqsubseteq (q^{1/2},1]$, and $\epsilon$ takes the values 
$1$ and $-1$. 
Rep\-re\-sen\-ta\-tions 
labeled by different sets of parameters (within unitary equivalence)
or belonging to 
different series are not unitarily equivalent. 
A representation of this list is 
irreducible if and only if $\K=\C$. In this case, $A$, $B$ and  $H$ 
are real numbers such that $A\in (q,1]$, $B\in (q,1]$,
and $H\in (q^{1/2},1]$, respectively. 
%
%
\subsection{Heisenberg representations 
of the cross product algebras {\mathversion{bold}$
 \cU_q(\mathrm{e}_2)\!\ltimes\!\cO (\mathrm{E}_q(2))$} and 
{\mathversion{bold}$\cU_q(\mathrm{e}_2)\!\ltimes\!\cO (\C_q)$} }
                                                     \label{3.5}
%
Let $\dC [u,v]$ be the algebra of complex Laurent 
polynomials $p(u,v)=\sum \alpha_{nk} u^n v^k$ 
in two commuting variables $u,v$ and let $\cF(\dR^+)$ 
be the algebra of all complex Borel functions $f=f(r)$ 
on $\dR^+=(0,+\infty)$ such that $f$ is locally bounded, that is, 
the restriction of $f$ to any compact subset 
contained in $\dR^+$ is bounded.
We denote by  $\cF (\mathrm{E}_q(2))$ the $\ast$-al\-ge\-bra generated 
by the two algebras $\dC [u,v]$ and $\cF (\dR^+)$ 
with cross commutation relations and involution 
\begin{equation}                                         \label{35}
u^j v^k f(r)=f(q^{-k}r) u^j v^k,\ \ 
(u^j v^k f(r))^\ast = \bar{f} (r) v^{-k} u^{-j},\ \ 
j,k\in \dZ,\ \, f\in  \cF (\dR^+).
\end{equation}

We introduce  a right action $\anf$ of $\cU_q(\mathrm{e}_2)$ 
on $\cF (\mathrm{E}_q(2))$ by 
\begin{align}\label{Krone}
&u^j v^k f(r) \anf E 
= q^{(j-k-3)/2} \lambda^{-1} 
u^{j+1} v^{k+1} (f(r)-q^{-j} f (qr))r^{-1},\\
\label{Zweig}
&u^j v^k f(r) \anf F = q^{(-j-k+1)/2} \lambda^{-1} 
u^{j-1} v^{k-1} (q^j f(r)-f(q^{-1}r))r^{-1},\\
\label{Ast}
&u^j v^k f(r) \anf K = q^{(j+k)/2} u^{j} v^{k} f(r), 
\end{align}
where  $j,k\in\dZ$ and $f\in\cF (\dR^+)$. 
Straightforward computations 
show that these formulas define indeed 
a right action of $\cU_q (\mathrm{e}_2)$ 
such that $\cF(\mathrm{E}_q(2))$ becomes 
a right $\cU_q(\mathrm{e}_2)$-mod\-ule $\ast$-al\-ge\-bra. 

From \rf[35], it follows that there is 
a $\ast$-iso\-mor\-phism $\phi$ 
from  $\cO (\mathrm{E}_q(2))$ onto the $\ast$-sub\-al\-ge\-bra 
of $\cF( \mathrm{E}_q(2))$ generated by $v$ and $ur$ 
such that $\phi(v)=v$ and $\phi (n)=ur$. 
From (\ref{Krone})--(\ref{Ast}), we conclude that 
$\phi$ intertwines the $\cU_q(\mathrm{e}_2)$-ac\-tion, that is, 
$x\anf f=\phi (x)\anf f$ for all $f\in \cU_q(\mathrm{e}_2)$ and 
$x \in \cO(\mathrm{E}_q(2))$. 
Identifying  $x\in\cO (\mathrm{E}_q(2))$ 
with  $\phi (x) \in \cF(\mathrm{E}_q(2))$, we can consider 
$\cO(\mathrm{E}_q(2))$ as 
a right $\cU_q(\mathrm{e}_2)$-mod\-ule 
$\ast$-sub\-al\-ge\-bra of $\cF(\mathrm{E}_q(2))$,  
and $\cU_q(\mathrm{e}_2) \lti \cO (\mathrm{E}_q(2))$ 
becomes a $\ast$-sub\-al\-ge\-bra of the corresponding 
cross product algebra $\cU_q(\mathrm{e}_2)\lti \cF (\mathrm{E}_q(2))$.

Let us briefly indicate how Equations (\ref{Krone})--(\ref{Ast}) 
have been derived. 
Using the definition of the  right actions of the generators $E$, $F$, $K$ 
on $n$, $n^\ast$ and  the fact that $\cO (\mathrm{E}_q(2))$ 
is a $\cU_q(\mathrm{e}_2)$-mod\-ule algebra, 
one obtains (\ref{Krone})--(\ref{Ast}) by straightforward calculations  
if  $j=k=0$ and $f(r)$ is a 
polynomial in $r^2 := n^\ast n$. Now 
we {\it postulate} that the formulas 
hold for arbitrary functions  $f(r)$.
Setting $u=nr^{-1}$ and $u^{-1}=r^{-1}n^\ast$, one can compute 
the actions of  $E$, $F$, $K$ 
on arbitrary elements $u^j v^k f(r)$ obtaining 
Equations (\ref{Krone})--(\ref{Ast}).  
On the other hand, 
since $n^j=u^j r^j$, $n^{\ast j}=u^{-j} r^j$, 
we rediscover the actions 
of the generators $E,F,K$ on the vector 
space basis $\{n^j v^k, n^{\ast l} v^k\,;\, j\in \dN_0,\  
l\in \dN,\  k\in \dZ\}$ of $\cO (\mathrm{E}_q(2))$
from Equations (\ref{Krone})--(\ref{Ast}). 

We turn now to the description of
a $\cU_q(\mathrm{e}_2)$-in\-var\-i\-ant 
positive linear functional $h_{\mu_0}$. Let $\mu_0$ 
be a finite positive Borel measure 
on the interval $(q,1]$. We extend $\mu_0$ to a positive 
Borel measure on $\dR^+$, denoted by $\mu$, 
such that 
$\mu (q^k \cM)=q^k\mu_0 (\cM)$ for  $k\in\dZ$, $\cM \subseteq (q,1]$. 
Let $\cF_0(\mathrm{E}_q(2))$ denote the $\ast$-sub\-al\-ge\-bra 
of $\cF(\mathrm{E}_q(2))$ generated by all elements $p(u,v) f(r)$, 
where $p\in\dC [u,v]$ and  $f\in\cF(\dR^+)$ has compact support.  
Define
\begin{equation}\label{Blatt}
h_{\mu_0} (p(u,v) f(r)) 
=\into_{\dT^2}p(u,v)\dd u\dd v{\into^{+\infty}_0} f(r)r\,\dd\mu (r).
\end{equation}
Using the above formulas (\ref{Krone})--(\ref{Ast}) 
for the actions of $E,F$ and $K$, one easily veri\-fies 
that $h_{\mu_0}(x\anf Z)=\varepsilon (Z)\hs h_{\mu_0} (x)$ 
for $x=u^n v^k f(r)$ and $Z=E,F,K,K^{-1}$. 
Since $\cF_0 (\mathrm{E}_q(2))$ is 
a $\cU_q (\mathrm{e}_2)$-mod\-ule algebra, 
this implies 
$h_{\mu_0} (x\anf Z)
=\varepsilon (Z) h_{\mu_0} (x)$ 
for all $x\in\cF_0 (\mathrm{E}_q(2))$ 
and $Z\in\cU_q(\mathrm{e}_2)$. 
Carrying out the integration over $\dT^2$, we obtain 
for $x=\sum_{j,k} \alpha_{jk} u^j v^k f_{jk} (r)$ 
and $y=\sum_{j,k} \beta_{jk} u^j r^k g_{jk} (r)$ 
from  $\cF_0 (\mathrm{E}_q(2))$
\begin{equation}\label{Stamm}
h_{\mu_0} (y^\ast x)
=\sumop_{j,k} \alpha_{jk}\overline{\beta_{jk}} 
{\into^{+\infty}_0} f_{jk} (r) \overline{g_{jk}} (r)r \,\dd\mu(r)
\end{equation}
This shows that  $h_{\mu_0}(x^\ast x)\ge 0$ 
for all $x\in \cF_0 (\mathrm{E}_q(2))$. 
Thus  $h_{\mu_0}$ is a $\cU_q(\mathrm{e}_2)$-in\-var\-i\-ant 
positive linear functional on the right 
$\cU_q(\mathrm{e}_2)$-mod\-ule 
$\ast$-al\-ge\-bra $\cF_0 (\mathrm{E}_q(2))$. 

Finally, we describe the Heisenberg representation of the cross product 
algebra $\ue\lti\cF (\mathrm{E}_q(2))$ associated with the 
invariant positive linear functional $h\equiv h_{\mu_0}$. 
From (\ref{Stamm}), 
it follows that the underlying Hilbert space is the tensor 
product of the Hilbert spaces $\cL^2 (\dT^2)$ and $\cL^2 (\dR^+, r \dd\mu )$. 
The actions of the generators $n, v\hsp\in\hsp \cO (\mathrm{E}_q(2))$, 
$f(r)\hsp \in\hsp \cF (\mathrm{E}_q(2))$ and $Z\hsp \in\hsp \ue$ 
on $\cL^2 (\dT^2)\hsp \otimes\hsp \cL^2 (\dR^+, r \dd\mu )$
are given by 
\begin{align*}
 &\pi_h (n) (u^jv^k \zeta (r)) = q^k u^{j+1}v^k r \zeta (r),\quad
 \pi_h (v)( u^jv^k \zeta (r)) = u^jv^{k+1}   \zeta (r),\\
&\pi_h (f(r))( u^jv^k \zeta (r))\hsp =\hsp  u^jv^k f(q^k r)\zeta (r), \ \ 
\pi_h (Z)( u^jv^k \zeta (r))\hsp  =\hsp ( u^jv^k \zeta (r))\anf S^{-1}(Z), 
\end{align*}
where $j,k\in\Z$ and $\zeta (r)\in \cF_0 (\mathrm{E}_q(2))$.

Let $\Hh={\oplus^\infty_{m,k,l=-\infty}}\Hh_{mkl}$, where each 
Hilbert space $\Hh_{mkl}$ is  $\cL^2 ((q,1], r\dd\mu_0)$. 
Define a linear operator 
$W:\Hh\rightarrow \cL^2 (\dT^2)\otimes \cL^2 (\dR^+, r\dd\mu)$ by 
$$
   W\zeta_{mkl}:=q^{m} u^{k}v^{l}\zeta(q^{m}r),\quad 
\zeta\in \cL^2 ((q,1], r\dd\mu_0),\ \, m,k,l\in\Z,
$$
From 
$$
{\into^{q^{-m}}_{q^{-m+1}}}|q^{m}\zeta(q^mr)|^2 r\,\dd\mu (r)
=\hspace{-1.5pt}
{\into^{q^{-m}}_{q^{-m+1}}}|\zeta(q^mr)|^2 q^{m}r\,\dd\mu (q^{m}r)
=\hspace{-1.5pt}{\into^{1}_{q}}|\zeta(r)|^2 r\,\dd\mu (r),
$$
it follows that $W$ is isometric. Since 
$\Lin \{ W\zeta_{mkl};
\zeta\hs{\in}\hs \cL^2 ((q,1], r\dd\mu_0),\, m,k,l\hs{\in}\hs\Z\}$ 
is dense in $\cL^2 (\dT^2)\otimes \cL^2 (\dR^+, r\dd\mu)$, we conclude 
that $W$ is  a unitary operator. 
Hence the Heisenberg representation on 
$\cL^2 (\dT^2)\otimes \cL^2 (\dR^+, r\dd\mu)$ is unitarily equivalent to 
a $\ast$-re\-pre\-sen\-ta\-tion on $\Hh$. Straightforward calculations show 
that the action of $\cU_q(\mathrm{e}_2)\lti\cO (\mathrm{E}_q(2))$ 
is determined by the following formulas: 
\begin{align*}
& v\zeta_{mkl}=\zeta_{m,k,l+1},\quad  
      n \zeta_{mkl} =q^{-m+l}Q \zeta_{m,k+1,l},\quad  
K\zeta_{mkl}=q^{-(k+l)/2} \zeta_{mkl},\\
&   F\zeta_{mkl} = -q^{m+(k-l-1)/2}\lambda^{-1}Q^{-1}\zeta_{m,k-1,l-1}
+q^{m-(k+l+1)/2}\lambda^{-1}Q^{-1}\zeta_{m-1,k-1,l-1},\\
&   E\zeta_{mkl} = -q^{m+(k-l-1)/2}\lambda^{-1}Q^{-1}\zeta_{m,k+1,l+1}
+q^{m-(k+l+1)/2}\lambda^{-1}Q^{-1}\zeta_{m+1,k+1,l+1},
\end{align*}
where $Q$ denotes the multiplication operator on the Hilbert space  
$\cL^2 ((q,1], r\dd\mu_0)$,
that is, $Q\zeta(r):=r\zeta(r)$. 
We carry out the computations only for the generator $E$, the other 
formulas are proved analogously. Recall from Subsection \ref{basics-Heis} that 
$\pi_h(E) (p(u,v) f(r))= (p(u,v) f(r))\anf S^{-1}(E)$. 
This gives 
\begin{align*}
E\zeta_{mkl}&=W^{-1}(-q^{m+1} (u^{k}v^{l}\zeta(q^mr))\anf E)\\
&=W^{-1}\big(-q^{m+(k-l-1)/2}\lambda^{-1}u^{k+1}v^{l+1}
(r^{-1}\zeta(q^mr)-q^{-k}r^{-1}\zeta(q^{m+1}r)) \big)\\
&=-q^{m+(k-l-1)/2}\lambda^{-1}Q^{-1}\zeta_{m,k+1,l+1}
+q^{m-(k+l+1)/2}\lambda^{-1}Q^{-1}\zeta_{m+1,k+1,l+1}
\end{align*}
as asserted. Now let $\beta\in\Z$ such that 
$-q^{-\beta-1}\lambda^{-1}\in(q,1]$. 
We define a unitary transformation $U$ on $\Hh$ by 
renaming the indices in the following way: 
$$
   U\eta_{mkl}:=(-1)^m\zeta_{\beta+k+l,\beta-m+k,\beta-m+k+l},\ \ 
   U^{-1}\zeta_{mkl}=(-1)^m\eta_{m-l,-\beta+m+k-l,-k+l}.
$$
Computing $U^{-1}fU\eta_{mkl}$ for the generators $f=v,n,K,E,F$ 
shows the  representation 
of\ \,$\U_q(\mathrm{e}_2)\lti\cO(\mathrm{E}_q(2))$
on $\Hh$ is 
unitarily equivalent to the $\ast$-re\-pre\-sen\-ta\-tion 
$(II)_{Q,-q^{-\beta-1}\lambda^{-1}I,I,1}$
from the preceding subsection.

Our next aim is the construction of 
a $\cU_q(\mathrm{e}_2)$-in\-var\-i\-ant 
positive linear functional
for the quantum complex plane $\C_q$. 
We proceed in a similar manner as in the case of $\mathrm{E}_q(2)$. 
Let $\dC[w]$ denote the algebra of complex Laurent polynomials in $w$ 
and let $\cF (\C_q)$ be the $\ast$-al\-ge\-bra 
generated by the two algebras $\cF (\dR^+)$ and 
$\dC[w]$ with cross 
commutation relation  $w^k f(r) = f(q^{-k} r) w^k$ 
and involution $(w^k f(r))^\ast = \bar{f} (r) w^{-k}$,  
where $k\in\dZ$ and $f\in\cF (\dR^+)$. 
We turn $\cF(\C_q)$ into a right $\cU_q(\mathrm{e}_2)$-mod\-ule 
$\ast$-algebra with right $\cU_q(\mathrm{e}_2)$-ac\-tion $\anf$ by setting  
\begin{align}\label{Krone1}
&w^k f(r) \anf E 
= q^{-3/2} \lambda^{-1} 
w^{k+1} (f(r)-q^{-k} f(qr))r^{-1},\\
\label{Zweig1}
&w^k f(r) \anf F = q^{1/2} \lambda^{-1} 
 w^{k-1} (f(r)-q^{-k} f(q^{-1}r))r^{-1},\\
\label{Ast1}
&w^k f(r) \anf K = q^{k} w^{k} f(r), 
\end{align}
where  $k\in\dZ$ and $f\in\cF (\dR^+)$. 
There is a $\ast$-iso\-mor\-phism $\phi$ 
from  $\cO (\C_q)$ onto the $\ast$-sub\-al\-ge\-bra 
of $\cF( \C_q)$ generated by $wr$ such that  $\phi(z)=wr$
which intertwines the $\cU_q(\mathrm{e}_2)$-ac\-tion.
Again, let  $\cF_0(\C_q)$ denote the $\ast$-sub\-al\-ge\-bra  of 
$\cF(\C_q)$ which is generated by $w$ and all functions 
$f\in \cF(\dR^+)$ with compact support. 
Recall that $\cO(\C_q)$ is a 
right  $\ue$-mod\-ule $\ast$-sub\-al\-ge\-bra of 
$\cO(\mathrm{E}_q(2))$ by identifying $z$ with $vn$. Comparing 
\rf[35] with the defining relations of $\cF(\C_q)$ and 
\rf[Krone]--\rf[Ast] with \rf[Krone1]--\rf[Ast1] shows that we 
can consider $\cF(\C_q)$ as a 
right  $\ue$-mod\-ule $\ast$-sub\-al\-ge\-bra of 
$\cF(\mathrm{E}_q(2))$ by identifying $w$ with $uv$.
From this, we deduce that the linear functional 
$\hat h_{\mu_0} $ defined by 
\begin{equation*}   
\hat h_{\mu_0} (p(w) f(r)) 
=\into_{\dT}p(w)\dd w{\into^{+\infty}_0} f(r)r\,\dd\mu (r)
\end{equation*}
is a $\cU_q(\mathrm{e}_2)$-in\-var\-i\-ant 
positive linear functional on $\cF_0(\C_q)$, where $\mu_0$ and $\mu$ 
are given as above. 
Moreover, the Heisenberg representation of $\ue\lti\cF(\C_q)$ 
associated with $\hat h_{\mu_0}$
is unitarily equivalent to the restriction of the Heisenberg 
representation of $\ue\lti\cF(\mathrm{E}_q(2))$ 
associated with $h_{\mu_0}$ to $\ue\lti\cF(\C_q)$. 
In particular, the Heisenberg representation of $\ue\lti\cO(\C_q)$ 
is unitarily equivalent to the representation 
on the subspace $\hat \Hh:={\oplus^\infty_{m,k=-\infty}}\Hh_{mkk}$ 
of $\Hh$, where the actions of $E$, $F$ and $K$ on 
$\zeta_{mkk}\in\Hh_{mkk}$ are given by above formulas and the 
action of $z$ is determined by 
$z \zeta_{mkk}=q^{-m+k}\zeta_{m,k+1,k+1}$. 
Setting $\eta_{nk}=(-1)^m\zeta_{\beta+k,\beta-m+k,\beta-m+k}$ and computing 
the actions of the generators $z$, $E$, $F$ and $K$ on $\eta_{nk}$, 
we obtain the formulas of the $\ast$-re\-pre\-sen\-ta\-tion 
$(II)_{Q,-q^{-\beta-1}\lambda^{-1}I,I}$
from Subsection \ref{rep-Cq}, where $\beta$ is defined as before. 

We summarize the main results of this subsection in the next proposition. 
\begin{thp} 
The Heisenberg representation 
of\ \,$\U_q(\mathrm{e}_2)\lti\cO(\mathrm{E}_q(2))$
associated with  $h_{\mu_0}$ is 
unitarily equivalent to the $\ast$-re\-pre\-sen\-ta\-tion 
$(II)_{Q,-q^{-\beta-1}\lambda^{-1}I,I,1}$
from Subsection \ref{rep-Eq2}
and the  Heisenberg representation 
of\ \,$\U_q(\mathrm{e}_2)\lti\cO(\C_q)$
associated with  $\hat h_{\mu_0}$ is 
unitarily equivalent to the $\ast$-re\-pre\-sen\-ta\-tion 
$(II)_{Q,-q^{-\beta-1}\lambda^{-1}I,I}$
from Subsection \ref{rep-Cq}, 
where $Q$ denotes the multiplication operator on the Hilbert space  
$\cL^2 ((q,1], r\dd\mu_0)$. 
\end{thp}

%
%
\section{Cross product algebras related to 
the quantum group {\mathversion{bold}$\mathrm{SU}_q(1,1)$}}
                                                     \label{sec-SU11}
%
\subsection{Definitions and ``decoupling'' of cross product algebras}
                                                       \label{SU-def}
For a moment, let $q$ be a complex number such that $q\ne 0,\pm 1$. 
First we repeat the definitions of the left and right cross product algebras 
of the Hopf algebra $\U_q(\mathrm{sl}_2)$ with 
the coordinate Hopf algebra $\cO(\mathrm{SL}_q(2))$  from [SW]. 
Recall that the algebra $\SUU$ has generators 
$a$, $b$, $c$, $d$ with defining relations 
\begin{equation}\label{orel}
ab=qba,\ \, ac=qca,\ \, bd=q db,\ \, cd=qdc,\ \,bc=cb,\ \, 
ad - qbc=da - q^{-1} bc=1,
\end{equation}
The Hopf algebra $\U_q(\mathrm{sl}_2)$ is 
generated by $E$, $F$, $K$, $K^{-1}$ 
with relations 
\begin{equation}\label{urel}
KK^{-1}\hsp\hspace{-0.04pt}=\hspace{-0.04pt}\hsp 
K^{-1}K\hspace{-0.04pt}\hsp=\hsp\hspace{-0.04pt} 1,\ KE\hsp=\hsp qEK,\  
KF\hsp=\hsp q^{-1} FK,\  
EF-FE\hsp=\hsp\lambda^{-1}(K^2-K^{-2}),
\end{equation}
and Hopf algebra structure 
\begin{align*}
&\Delta(K)=K\otimes K,\ \,
\Delta(E)=E\otimes K+K^{-1}\otimes E,\ \,
\Delta(F)=F\otimes K+K^{-1}\otimes F, \\
&\varepsilon (K)=1,\ \, \varepsilon (E) = \varepsilon (F)= 0,\ \, 
S(K)=K^{-1}, \ \, S(E)=-qE,\ \, S(F)=-q^{-1}F. 
\end{align*}

There is a dual pairing of Hopf algebras 
$\langle\cdot,\cdot\rangle : 
\U_q (\mathrm{sl}_2)\times \cO(\mathrm{SL}_q(2))
\rightarrow \C$ given on generators by 
$$                                        
\langle K^{\pm },d\rangle =\langle K^{\mp 1},a\rangle=q^{\pm 1/2},  
\quad \langle E,c\rangle =\langle F,b\rangle=1
$$
and zero otherwise. Using Equations (\ref{cross2}) and (\ref{cross3}), 
one derives the following cross commutation relations 
in the corresponding left and right cross product algebras.
\begin{align*}   
\U_q (\mathrm{sl}_2)\lti\cO(&\mathrm{SL}_q(2)):   \\
aK &=q^{-1/2} K a, & aE&=q^{-1/2} E a, & aF &=q^{-1/2} F a+K^{-1} c,\\
bK  &=q^{-1/2} Kb, &b E &=q^{-1/2} E b, & b F &=q^{-1/2} F b+K^{-1}d,\\
c K &= q^{1/2} K c, & cE&= q^{1/2} E c+K^{-1} a, & 
cF&=q^{1/2} F c,\\
dK&= q^{1/2} K d, &  d E &=q^{1/2} E d+K^{-1} b, &
d F&=q^{1/2} F d.\\[-24pt]
\end{align*}
\begin{align*}   
\cO(\mathrm{SL}_q(2))&\rti\U_q (\mathrm{sl}_2):\\
Ka&=q^{-1/2} aK, & 
Ea&=q^{1/2} aE+ bK, &  Fa&=q^{1/2} aF,\\
Kb&=q^{1/2} bK, &   E b&=q^{-1/2}b E,
 &  Fb&=q^{-1/2}b F+aK,\\
Kc &=q^{-1/2} cK, &   
Ec&=q^{1/2} cE+d K, &  Fc&=q^{1/2} cF,\\
Kd&= q^{1/2} dK, &  Ed&=q^{-1/2} E d,
 &  Fd&=q^{-1/2} dF+cK.
\end{align*}

Next we want to embed $\cO(\mathrm{SL}_q(2))$ 
into a larger algebra where $b$ and $c$ are invertible. 
For this reason, we consider the localization of 
$\cO(\mathrm{SL}_q(2))$ at the set  $\cS:=\{b^nc^m;n,m\in\N_0\}$. 
Clearly, $\cS$ is 
a left and right Ore set of $\cO(\mathrm{SL}_q(2))$, 
and the algebra $\cO(\mathrm{SL}_q(2))$ has no zero divisors. 
Therefore 
the localization $\hat{\cO}(\mathrm{SL}_q(2))$ 
of $\cO(\mathrm{SL}_q(2))$ at $\cS$ exists and
contains $\cO(\mathrm{SL}_q(2))$ as a subalgebra. 
Note that all elements $b^{-n} c^n$, $n\in\Z$, 
belong to the center of  $\hat{\cO} (\mathrm{SL}_q(2))$. 
From Theorem 3.4.1 in \cite{LR} or from Theorem 1.2 in \cite{Ko}, 
it follows that the algebra $\hat{\cO}(\mathrm{SL}_q(2))$ 
is a right (resp.\ left) $\U_q(\mathrm{sl}_2)$-mod\-ule algebra  
which contains $\cO(\mathrm{SL}_q(2))$ as 
a $\U_q(\mathrm{sl}_2)$-mod\-ule subalgebra. 
The actions of generators of $\U_q (\mathrm{sl}_2)$ 
on inverses  $b^{-1},c^{-1}$ are given by
\begin{align*}
b^{-1}\anf E&=0, & b^{-1} \anf F&=-q^{-1} db^{-2}, &
b^{-1}\anf K&=q^{1/2}b^{-1},\\
c^{-1}\anf E&=-qac^{-2},& c^{-1} \anf F&=0, & c^{-1} \anf K&=q^{-1/2}c^{-1},
\\[-24pt]
\end{align*}
\begin{align*}
E\ang b^{-1}&=0,& F\ang b^{-1}&=-q ab^{-2},&K\ang b^{-1}&= q^{1/2}b^{-1},\\
E\ang c^{-1}& =-q^{-1} dc^{-2},& F\ang c^{-1}&=0,& 
K\ang c^{-1}&=q^{1/2} c^{-1}. 
\end{align*}
Taking these formulas as definitions, 
one may also verify directly 
that $\hat{\cO}(\mathrm{SL}_q(2))$ is 
a $\U_q(\mathrm{sl}_2)$-mod\-ule algebra 
without using the results from \cite{LR}, \cite{Ko}.

Since $\hat{\cO} (\mathrm{SL}_q(2))$ is a left and right 
$\U_q(\mathrm{sl}_2)$-mod\-ule algebra, the cross product algebras  
$\U_q(\mathrm{sl}_2)\lti \hat{\cO} (\mathrm{SL}_q(2))$ 
and $\hat{\cO}(\mathrm{SL}_q(2))\rti\cU_q(\mathrm{sl}_2)$ 
are well defined. 
We can regard $\U_q(\mathrm{sl}_2)\lti \hat{\cO} (\mathrm{SL}_q(2))$ 
and $\hat{\cO}(\mathrm{SL}_q(2))\rti\cU_q(\mathrm{sl}_2)$ 
as algebras generated by $a$, $b$, $b^{-1}$, $c$, $c^{-1}$, $d$ 
and $E$, $F$, $K$, $K^{-1}$ with the defining relations 
of $\U_q(\mathrm{sl}_2)\lti \cO(\mathrm{SL}_q(2))$ 
and $\cO(\mathrm{SL}_q(2))\rti\U_q(\mathrm{sl}_2)$, respectively, 
and the additional relations 
\begin{equation}\label{bcrel}
bb^{-1}=b^{-1}b=cc^{-1}=c^{-1}c=1.
\end{equation}

The left cross product algebra is isomorphic to its right-handed 
counterpart. An isomorphism $\theta$ is given by 
$\theta(a)=a$, $\theta(b)=-qc$, $\theta(c)=-q^{-1}b$, $\theta(d)=d$, 
$\theta(E)=F$, $\theta(F)=E$ and $\theta (K)=K^{-1}$.

Next we show that, by passing to different generators, 
the defining relations of the cross product algebras 
simplify remarkably. The following lemma is proved by direct calculations.  
We restrict ourselves to the right-handed version. 
\begin{thl}
Set 
\begin{equation}\label{qrdef}
Q:=-q^{1/2}\lambda K^{-1} E - K^{-2} c^{-1} a,\quad 
R:= q^{1/2} \lambda FK^{-1}  - q db^{-1} K^{-2}.
\end{equation}
Then
\begin{align}\label{qrrel}
&xQ=Qx,\quad   xR=Rx,\quad  x\in\SUU, \\
\label{qrzrel}
&KQ=qQK,\quad KR=q^{-1} RK,\quad QR-q^2 RQ=1-q^2.
\end{align}
\end{thl}

By (\ref{qrdef}), we can write
\begin{equation}\label{efdef}
E=-q^{-1/2} \lambda^{-1} (KQ+K^{-1} c^{-1} a),\quad  
F=q^{1/2} \lambda^{-1} (RK+qd b^{-1} K^{-1}). 
\end{equation}
It is straightforward to check that 
$a$, $b$, $b^{-1}$, $c$, $c^{-1}$, $d$ and $Q$, $R$, $K$, $K^{-1}$
are generators of the cross product 
algebra\ \,$\U_q(\mathrm{sl}_2)\lti \hat{\cO}(\mathrm{SL}_q(2))$ 
satisfying the 
defining relations (\ref{orel}), (\ref{bcrel}), (\ref{qrrel}), 
(\ref{qrzrel}) and 
\begin{equation}\label{krel}
aK = q^{-1/2} Ka,\quad bK=q^{-1/2} Kb,\quad  cK=q^{1/2} Kc, \quad 
dK=q^{1/2} Kd.
\end{equation}
Let $\U$ denote the subalgebra of 
$\U_q(\mathrm{sl}_2)\lti \hat{\cO}(\mathrm{SL}_q(2))$ 
generated by $Q$, $R$, $K$, $K^{-1}$. 
Then 
$\U_q(\mathrm{sl}_2)\lti\hat{\cO}(\mathrm{SL}_q(2))$ 
can be considered as the algebra generated by the subalgebras 
$\hat{\cO}(\mathrm{SL}_q(2))$ and\ \,$\U$ 
with  ``almost decoupled''
cross relations (\ref{qrrel}) and (\ref{krel}).

The main advantage of the new generators $Q$ 
and $R$ is that they commute with the elements 
of the algebra $\hat{\cO}(\mathrm{SL}_q(2))$ by (\ref{qrrel}). 
This fact and the form of these generators can also 
be derived from Lemma \ref{2.2} applied 
to the right coideals $\V=\Lin \{EK,\varepsilon\}$,   
$\V' = \Lin \{ FK,\varepsilon\}$ and the set of generators 
$\X_0=\{a,b,c,d\}$.
Indeed, for $v=EK$, condition (\ref{xvcon}) means that 
\begin{align*}
 &a\rho (EK)=q^{-1} \rho (EK) a, & b\rho (EK)&=q^{-1} \rho (EK) b,\\
&c\rho (EK)= q\rho (EK) c+q^{-1/2} a, & d\rho (EK)&=q\rho (EK) d+q^{-1/2} b.
\end{align*}
Therefore, setting $\rho (EK)=-q^{-3/2} \lambda^{-1} c^{-1}a$ 
and $\rho (\varepsilon)=1,$ (\ref{xvcon}) 
is satisfied for $\V=\Lin\{ EK,\varepsilon\}$, and we get
\begin{align*}
\xi (EK) &= \rho (EK) S(K^2)  +  \rho (\varepsilon) S(EK)
=-q^{-3/2} \lambda^{-1} c^{-1} aK^{-2}  - q K^{-1} E \\
&= q^{1/2} \lambda^{-1} Q.
\end{align*}
Similarly, with $\rho (FK)=q^{1/2} \lambda^{-1} db^{-1}$, 
we obtain $\xi (FK) =-q^{-1/2} \lambda^{-1} R$. 

It might be worth to mention that 
there is no algebra homomorphism 
$\varphi$ from $\U_q(\mathrm{sl}_2)\lti \hat{\cO} (\mathrm{SL}_q(2))$ 
into $\hat{\cO}(\mathrm{SL}_q(2))$ such that $\varphi(x)=x$ 
for all $x\in\hat{\cO}(\mathrm{SL}_q(2))$. Indeed, one can easily 
show that there is no element 
$\varphi (K^2)\in\hat{\cO}(\mathrm{SL}_q(2))$ 
such that $b\varphi (K^2)=q^{-1} \varphi(K^2)b$ 
and  $c\varphi(K^2)=q\varphi(K^2)c$. 
Hence the results of \cite{F} do not apply to the 
cross product algebra 
$\cU_q(\mathrm{sl}_2)\lti\hat{\cO}(\mathrm{SL}_q(2))$.

Let us remark that 
the generators $Q$ and $R$ behave nicely  
under the involutions of the three real 
forms of $\U_q(\mathrm{sl}_2))$ and ${\cO}(\mathrm{SL}_q(2))$. 
For $q$ real, we have $Q^\ast=-R$ and $Q^\ast=R$ 
in the $\ast$-al\-ge\-bras 
$\U_q(\mathrm{su}_2)\lti \cO(\mathrm{SU}_q(2))$ and 
$\suu\lti\SUU$, respectively. The third relation of (\ref{qrzrel}) 
reads then  $QQ^\ast -q^2 Q^\ast Q=\pm (1-q^2)$ 
with the minus sign in the first case and plus in the second. 
The representations of this relation 
are described in Lemma \ref{2.2}.
For $|q|=1$, we have  
$(q^{1/2}Q)^\ast=q^{1/2} Q$ and  $(q^{-1/2} R)^\ast=q^{-1/2} R$ 
in the $\ast$-al\-ge\-bra 
$\U_q(\mathrm{sl}_2(\R))\lti\hat{\cO} (\mathrm{SL}_q(2,\R))$.  
Here the involutions of the Hopf $\ast$-al\-ge\-bras 
$\U_q(\mathrm{sl}_2(\R))$ and $\cO (\mathrm{SL}_q(2,\R))$ 
are defined by 
$E^\ast=-qE$, $F^\ast=-q^{-1}F$, $K^\ast =K$ and $a^\ast =a$, 
$b^\ast =b$, $c^\ast=c$, $d^\ast=d$ so that 
$\langle\cdot,\cdot\rangle$ is a dual pairing of 
Hopf $\ast$-al\-ge\-bras and 
$\U_q(\mathrm{sl}_2(\R))\lti\cO (\mathrm{SL}_q(2,\R))$ is indeed 
a $\ast$-al\-ge\-bra. 
In all three cases,
the algebra $\U$ generated by $Q$, $R$, $K$, $K^{-1}$ is a $\ast$-al\-ge\-bra.

From now, we suppose again that $q\in (0,1)$. 
We are interested in the real forms 
$\SUU$ and $\suu$. 
On generators, the involution is given by 
$$
a^\ast=d,\quad b^\ast=qc,\qquad  E^\ast=-F,\quad K^\ast=K.
$$
The pairing $\langle\cdot,\cdot\rangle$ defined above is a 
dual pairing of Hopf $\ast$-al\-ge\-bras, so  
the cross pro\-duct algebras 
$\U_q(\mathrm{su}_{1,1})\lti \hat{\cO}(\mathrm{SU}_q(1,1))$  
and $\hat{\cO}(\mathrm{SU}_q(1,1))\rti \U_q(\mathrm{su}_{1,1})$ 
are  $\ast$-al\-ge\-bras with involutions induced from 
the $\ast$-al\-ge\-bras $\cO(\mathrm{SU}_q(1,1))$ and 
$\U_q(\mathrm{su}_{1,1})$. 
The mapping $\theta$ realizing the isomorphism of the left and 
right cross product algebras is a $\ast$-iso\-mor\-phism.  
Hence the $\ast$-al\-ge\-bras  
$\U_q(\mathrm{su}_{1,1})\lti \hat{\cO}(\mathrm{SU}_q(1,1))$ 
and $\hat{\cO}(\mathrm{SU}_q(1,1))\rti \U_q(\mathrm{su}_{1,1})$  
are $\ast$-iso\-mor\-phic. 

Now we turn to cross product algebras related to the quantum disc. 
The quantum disc algebra $\qd$ is defined as the $\ast$-al\-ge\-bra generated 
by $z$ and $z^\ast $ with relation
\begin{equation}                                  \label{zzstar}
               z^\ast z-q^2zz^\ast =1-q^2.
\end{equation}

A left action $\ang$ which turns $\qd$ into a $\suu$-mod\-ule $\ast$-al\-ge\-bra 
appears in \cite{KL} and \cite{VK3}.  On generators, it takes the 
form 
\begin{align*}                                     
 K^{\pm1}\ang z&=q^{\mp 1}z, &  E\ang z&=q^{1/2}, &  
 F\ang z&=-q^{-1/2}z^2,\\
 K^{\pm1}\ang z^\ast &=q^{\pm1}z^\ast , &  
E\ang z^\ast &=-q^{1/2}z^{\ast 2}, &  
 F\ang z^\ast &=q^{-1/2}.
\end{align*}
A right action $\anf$ can be obtained from $\ang$ by applying 
Lemma \ref{r-l} with the algebra 
anti-automorphism and coalgebra homomorphism
$\phi:\suu\rightarrow\suu$ given by $\phi(K)=K$,
$\phi(E)=qF$, $\phi(F)=q^{-1}E$. 
From \rf[ract], we derive 
\begin{align*}                                     
 z\anf K^{\pm1}&=q^{\mp 1}z, &  z\anf E&=-q^{1/2}z^2, &  
 z\anf F&=q^{-1/2},\\
 z^\ast \anf K^{\pm1}&=q^{\pm 1}z^\ast , &  z^\ast \anf E&=q^{1/2}, &  
 z^\ast \anf F&=-q^{-1/2}z^{\ast 2}.
\end{align*}
These formulas lead to the following cross commutation relations in 
the corresponding cross product algebras. 
\begin{align*}   
\cU_q(\mathrm{su}&_{1,1})\lti\makebox[0cm][l]{$\qd:$}   \\
zK&=q^{-1}Kz, &  zE&=q^{-1}Ez-q^{1/2}K^{-1}z^2, &  
zF&=q^{-1}Fz+q^{-1/2}K^{-1},\\
z^\ast K&=qKz^\ast , &  z^\ast E&=qEz^\ast +q^{1/2}K^{-1},
 &  z^\ast F&=qFz^\ast -q^{-1/2}K^{-1}z^{\ast 2}.\\[-24pt]
\end{align*}
\begin{align*}
\cO(\mathrm{U}_q)\rti&\cU_q(\mathrm{su}_{1,1}):   \\
Kz&=q^{-1}zK, &  Ez&=qzE+q^{1/2}K, &  Fz&=qzF-q^{-1/2} z^2K,\\
Kz^\ast&=qz^\ast K, &  Ez^\ast&=q^{-1}z^\ast E-q^{1/2} z^{\ast 2} K, 
 &  Fz^\ast&=q^{-1}z^\ast F+q^{-1/2} K.       
\end{align*}
The $\ast$-al\-ge\-bras 
$\U_q(\mathrm{su}_{1,1})\lti\cO(\rmU_q)$ and 
$\cO(\rmU_q)\rti \U_q(\mathrm{su}_{1,1})$ are 
$\ast$-iso\-morphic with a $\ast$-iso\-mor\-phism $\psi$ 
determined by $\psi(z)=z$, $\psi(K)=K^{-1}$ and $\psi(E)=-F$.

There is also a ``decoupling''  for the cross commutation 
relations of the cross product algebra 
$\U_q(\mathrm{su}_{1,1})\lti\cO(\rmU_q)$. 
The elements
$$
S\hsp:=\hsp q^{1/2} \lambda FK^{-1}  - qz^\ast K^{-2}, \ \,
S^\ast\hsp=\hsp -q^{-1/2} \lambda K^{-1} E - q K^{-2} z, \ \,
T\hsp:=\hsp K^{-2} (1 - z^\ast z)
$$ 
of the algebra $\U_q(\mathrm{su}_{1,1})\lti \cO(\rmU_q)$ 
satisfy the relations 
\begin{align}\label{zstrel}
&zS\hsp =\hsp Sz,\  \,z^\ast S\hsp =\hsp Sz^\ast, \ \, 
zS^\ast\hsp =\hsp S^\ast z,\ \, z^\ast S^\ast\hsp =\hsp S^\ast z^\ast, \ \,
zT\hsp =\hsp Tz,\ \, z^\ast T\hsp =\hsp Tz^\ast \\
\label{strel}
&ST=q^{-2} TS, \quad  S^\ast T= q^2 TS^\ast, \quad 
S^\ast S- q^2 SS^\ast = 1-q^2.
\end{align}
Thus, the $\ast$-sub\-al\-ge\-bra generated by $S$, $S^\ast$ and  
$T\hsp =\hsp T^\ast$ commutes with the $\ast$-sub\-al\-ge\-bra $\cO(\rmU_q)$ 
of $\U_q(\mathrm{su}_{1,1})\lti\cO(\rmU_q)$.  
These two $\ast$-sub\-al\-ge\-bras generate 
a large part but not the whole of the 
cross product algebra $\U_q(\mathrm{su}_{1,1})\lti\cO(\rmU_q)$. 
An alternative set of generators 
of $\U_q(\mathrm{su}_{1,1})\lti\cO(\rmU_q)$ 
is $z$, $z^\ast$, $S$, $S^\ast$, $K$, $K^{-1}$. 
Then the defining relations are (\ref{zzstar}), 
the corresponding relations of (\ref{zstrel}) and (\ref{strel}), and 
$$
KK^{-1}\hsp =\hsp K^{-1}K\hsp =\hsp 1,\  Kz\hsp =\hsp qzK,\  
z^\ast K\hsp =\hsp qKz^\ast,\   
SK\hsp =\hsp qKS,\  
K S^\ast \hsp =\hsp qS^\ast K.
$$ 
The form of the generators $S$, $S^\ast$, $T$ and the 
fact that they commute with the algebra $\cO(\rmU_q)$ 
can also be obtained from Lemma \ref{2.2} applied 
to right coideals $\V=\Lin \{FK,\varepsilon\}$, 
$\V'=\Lin\{EK,\varepsilon\}$, $\V''=\Lin \{K^2\}$ and the set
$\X_0=\{z,z^\ast\}$. 
Equation (\ref{xvcon}) is satisfied if we set 
$\rho(FK)=q^{1/2} \lambda^{-1} z^\ast$, $\rho(EK)
=-q^{-1/2} \lambda^{-1} z$, $\rho(K^2)=1 - z^\ast z$  
so that  $\xi(FK)=-q^{-1/2} \lambda^{-1} S$, $\xi(EK)=q^{1/2} 
\lambda^{-1} S^\ast$ and $\xi(K^2)=T$. 

Analogously to Subsection \ref{sec-def}, one can also consider 
a cross product $\ast$-sub\-al\-ge\-bra 
$\U_0\lti\SUU$ of $\suu\lti\SUU$, where 
$\U_0\subset\suu$ is the unital $\ast$-al\-ge\-bra 
generated by the quantum tangent space of a 
left-covariant first order differential $\ast$-cal\-cu\-lus on $\SUU$. 
For Woronowicz's 3D-calculus, this cross product algebra 
and its representations have been studied 
in \cite{EW1}. 
%
\subsection{\!\!\!Representations 
of the {\mathversion{bold}$\ast$}-al\-ge\-bra 
{\mathversion{bold}$ \suu\!\ltimes\!\SUU$} }
                                                \label{sec-SUU}
By applying Lemma \ref{L3} to the relation 
$aa^\ast -q^2a^\ast a=1-q^2$, using the identity $c^\ast c=a^\ast a-1$, 
taking the polar decomposition of the closed operator $c$, 
and arguing as in Subsection \ref{sec-U0Cq}, 
one easily shows that 
any $\ast$-re\-pre\-sen\-tation of $\SUU$ 
is unitarily equivalent to a representation on the orthogonal sum 
$\G\oplus\Hh$ of Hilbert spaces $\G$ and $\Hh$ determined by 
\begin{align*}
& a=v,\quad d=v^\ast,
\quad  b=c=0\quad \mbox{on}\ \,\G,\\  
&a\eta_n=(1+q^{2n} A^2)^{1/2}\eta_{n-1},\quad 
   d\eta_n=(1+q^{2(n+1)} A^2)^{1/2}\eta_{n+1}, \\ 
&b\eta_n=q^{n+1} Aw^\ast \eta_n, \quad  c\eta_n=q^n Aw\eta_n
\quad \mbox{on}\ \,
\Hh=\displaystyle\mathop{\oplus}_{n=-\infty}^\infty \Hh_n,\ \,\Hh_n=\K.
\end{align*}
Here, $w$ is a unitary and $A$ is a self-adjoint operator on a 
Hilbert space $\K$ satisfying $wA=Aw$ and  
$\sigma(A)\sqsubseteq (q,1]$, and  $v$ is a unitary operator on $\G$
(see \cite{EW1}).

On $\ast$-rep\-re\-sen\-ta\-tions of $\suu\lti\SUU$, we impose the 
following regularity conditions. We assume that the restriction to $\SUU$ 
is of the form described above and that there exist 
dense linear subspaces $\E$ of $\G$ and $\D_0$ 
of $\Hh_0$ such that $v\E=\E$, $w\D_0=\D_0$, $A\D_0=\D_0$, and 
$\E\oplus\D$ is invariant under the actions of $a$, $b$, $c$, $d$,  
$E$, $F$, $K$ and $K^{-1}$, where 
$\D=\Lin\{\eta_n; \eta \in \D_0,~n\in\Z\}$.

First we show that $\G=\{0\}$. 
Note that $\G=\Ker b=\Ker c$.
From $KEb=qbKE $, it follows that the operator $KE$ leaves $\E$ 
invariant. Thus $v\eta=a\eta=(q^{-1/2}cKE-q^{1/2}KE c)\eta=0$ 
for all $\eta\in\E$. Since $v$ is unitary and $\E$ is dense in $\G$, 
we have $\G=\{0\}$. 

We assume that the commutation relations of $K$ with $E$, $F$  
and the generators of $\SUU$ hold in strong sense. Then 
it follows that the subspaces 
of $\Hh$ on which $K>0$ and $K<0$ are reducing. 
Therefore we may assume that $K= \epsilon |K|$ with  
$\epsilon\in\{1,-1\}$. 
From $c^\ast cK = K c^\ast c$, it follows that $K$ leaves 
each Hilbert space $\Hh_n$ invariant. 
Hence there exist positive self-adjoint 
operators $H_n$ on $\Hh_0$ such that 
the action of $K$ on 
$\Hh=\oplus_{n=-\infty}^\infty \Hh_n$, $\Hh_n=\Hh_0$, is given by 
$H\eta_n=\epsilon H_n\eta_n$, and each $H_n$ commutes strongly with $A$. 
The relation $a K= q^{-1/2}Ka$ applied 
to vectors $\eta_{n+1}\in\D\cap\Hh_{n+1}$
implies $H_{n+1}= q^{-1/2}H_{n}$ since 
$H_n\alpha_{n}(A)=\alpha_{n}(A)H_n$ and $\Ker\alpha_{n}(A)=\{0\}$.
Hence $H_n= q^{-n/2}H_0$. Moreover, from $cK=q^{1/2}Kc$, we derive 
$wH_0=q^{1/2}H_0w$. 

As $\G=\{0\}$, $b$ and $c$ are invertible. 
With $Q$ defined in \rf[qrdef], 
we assume that the relation $c^\ast c Q = Q c^\ast c$ 
holds in strong sense. Then it follows by arguments similar to those 
used in the preceding paragraph that $Q$ acts on $\Hh$ by 
$Q\eta_n=Q_0\eta_n$, where $Q_0$ denotes an operator 
on $\Hh_0$ satisfying $Q_0A=AQ_0$ and $wQ_0=Q_0w$. 
In addition, $KQ=qQK$ gives $H_0Q_0=qQ_0H_0$ and the last equation of 
\rf[qrzrel] yields $Q_0Q_0^\ast-q^2Q_0^\ast Q_0=(1-q^2)$. 

Summarizing, the $\ast$-rep\-re\-sen\-ta\-tions of $\suu\lti\SUU$ 
are obtained by solving 
the following operator equations on a dense linear subspace $\D_0$ of 
$\Hh_0$:
\begin{align}
&wA=Aw,\quad AQ_0=Q_0A,\quad wQ_0=Q_0w,                    \label{T1}\\ 
&AH_0=H_0A,\quad wH_0=q^{1/2}H_0w,\quad H_0Q_0=qQ_0H_0,   \label{T2}\\
&Q_0Q_0^\ast-q^2Q_0^\ast Q_0=(1-q^2),                     \label{T3}
\end{align}
where $w$ is a unitary, $H_0$ is a positive self-adjoint operator, and 
$A$ is a self-adjoint operator subjected to the spectral 
condition $\sigma(A)\sqsubseteq (q,1]$.  
In addition, we require 
that \rf[T2] and the first two relations of 
\rf[T1] hold in strong sense, and 
that $A\D_0=\D_0$ and $H_0\D_0=\D_0$. 

The rep\-re\-sen\-ta\-tions of \rf[T3] are listed 
in Lemma \ref{L3}. If $Q_0$ is given by the series $(I)$ or $(II)_B$, 
then $\Hh_0$ is a direct sum $\oplus_{k\in J}\Hh_{0k}$, 
$\Hh_{0k}=\Hh_{00}$, where $J=\N_0$ and $J=\Z$ for representations 
of type $(I)$ and $(II)_B$, respectively.   
A similar analysis 
as used to derive the identities \rf[T1] and \rf[T2] 
shows that the operators $w$, $A$ and $H$ act on 
$\Hh_0=\oplus_{k\in J}\Hh_{0k}$, $\Hh_{0k}=\Hh_{00}$, 
by $w\zeta_k=w_0\zeta_k$, $A\zeta_k=A_0\zeta_k$ and 
$H_0\zeta_k=q^{-k}H_{00}\zeta_k$, where $w_0$ is  unitary, 
$A_0$ is a self-adjoint operator satisfying $\sigma(A_0)\sqsubseteq (q,1]$ 
and $H_{00}$ is a positive self-adjoint operator on $\Hh_{00}$ such that 
$$
w_0A_0=A_0w_0,\quad A_0H_{00}=H_{00}A_0,\quad w_0H_{00}=q^{1/2}H_{00}w_0, 
$$
and, for representations of type $(II)_B$, 
$$
w_0B=Bw_0,\quad A_0B=BA_0,\quad H_{00}B=BH_{00}. 
$$
The rep\-re\-sen\-ta\-tions of these relations are described in 
Lemma \ref{1}. 
We choose the self-adjoint operator $B$ such that $1$ is not 
an eigenvalue, that is, $\sigma(B)\sqsubset[q,1)$.

Finally, suppose that $Q_0$ is given by the series $(III)_u$. 
Then $Q_0=u$ is a unitary operator on $\Hh_0$. 
Now we apply Lemma \ref{1} to the relations 
$wA=Aw$, $wH_0=q^{1/2}H_0w$ and $AH_0=H_0A$. 
It states that there exist commuting 
self-adjoint operators $A_0$ and $H_{00}$ on a Hilbert space $\Hh_{00}$
satisfying $\sigma(A_0)\sqsubseteq (q,1]$ and  
$\sigma(H_{00})\sqsubseteq (q^{1/2},1]$ 
such that $\Hh_0=\oplus_{k=-\infty}^{\infty}\Hh_{0k}$, 
$\Hh_{0k}=\Hh_{00}$, and the actions of $w$, $A$ and $H_0$ on $\Hh_0$
are determined by $w\zeta_k=\zeta_{k-1}$, $H_0\zeta_k=q^{k/2}H_{00}\zeta_k$ 
and $A\zeta_k=A_0\zeta_k$. From $w^2uH_0=Hw^2u$, it follows that $w^2u$ 
leaves each Hilbert space $\Hh_{0k}$ invariant. Hence we can write 
$w^2u\zeta_k=u_k\zeta_k$, where $u_k$ denotes a unitary operator 
on $\Hh_{0k}\,(=\Hh_{00})$. Accordingly, $u\zeta_k=u_k\zeta_{k+2}$. 
Since $wu=uw$, we have $u_k=u_{k-1}$, hence all $u_k$ are equal. 

Inserting the expressions derived for $a$, $b$, $c$, $d$, $Q$ and $R=Q^\ast$ 
into Equation \rf[qrdef], using the abbreviations $\alpha_n(t)$ 
and $\beta_n(t)$ 
introduced in the introduction, 
and renaming the operators, we obtain the following list of 
$\ast$-rep\-re\-sen\-ta\-tions of $\suu\lti\SUU$. 
\begin{align*}
(I.1&)_{A,H,\epsilon}: &
a\eta_{nkl}&=\alpha_n(A)\eta_{n-1,kl},\qquad  &
 d\eta_{nkl}&=\alpha_{n+1}(A)\eta_{n+1,kl}, \\
& &   b\eta_{nkl}&=q^{n+1}  A\eta_{nk,l+1},\qquad &  
c\eta_{nkl}&=q^n A\eta_{nk,l-1},\\
&F\eta_{nkl}=\makebox[0cm][l]{$ \lambda^{-1} q^{(-n-2k+l-1)/2}
                              \lambda_{k+1}\epsilon H\eta_{n,k+1,l} $}  \\
& & &+\makebox[0cm][l]{$ \lambda^{-1}q^{(n+2k-l+1)/2}\epsilon H^{-1}
                                      \beta_{n+1}(A)\eta_{n+1,k,l-1},$} \\
&E\eta_{nkl}=\makebox[0cm][l]{$ -\lambda^{-1}  q^{(-n-2k+l+1)/2}
                                  \lambda_k \epsilon H\eta_{n,k-1,l} $}\\
& & &-\makebox[0cm][l]{$ \lambda^{-1}q^{(n+2k-l-1)/2}\epsilon H^{-1}
                                        \beta_{n}(A) \eta_{n-1,k,l+1},$}\\
&K\eta_{nkl}=\makebox[0cm][l]{ 
$q^{(-n-2k+l)/2}\epsilon H\eta_{nkl}$ \quad  on \ \ $\Hh=
\displaystyle\mathop{\oplus}_{n,l=-\infty}^\infty 
\displaystyle\mathop{\oplus}_{k=0}^\infty 
\Hh_{nkl},\ \,\Hh_{nkl}=\K.$}\\[-24pt]
\end{align*}
\begin{align*}
(I.2&)_{A,B,H,\epsilon}: &
a\eta_{nkl}&=\alpha_n(A)\eta_{n-1,kl},\qquad  &
d\eta_{nkl}&=\alpha_{n+1}(A)\eta_{n+1,kl}, \\
& &   b\eta_{nkl}&=q^{n+1}  A\eta_{nk,l+1},\qquad &  
c\eta_{nkl}&=q^n A\eta_{nk,l-1},\\
&F\eta_{nkl}=\makebox[0cm][l]{$
 \lambda^{-1} q^{(-n-2k+l-1)/2}\alpha_{k+1}(B)\epsilon H\eta_{n,k+1,l} $}\\
& & &+\makebox[0cm][l]{$ \lambda^{-1}q^{(n+2k-l+1)/2}\epsilon H^{-1} 
                                        \beta_{n+1}(A)\eta_{n+1,k,l-1},$} \\
&E\eta_{nkl}=\makebox[0cm][l]{$ 
  -\lambda^{-1} q^{(-n-2k+l+1)/2}\alpha_k(B)  \epsilon H\eta_{n,k-1,l} $}\\
& & &-\makebox[0cm][l]{$ \lambda^{-1}q^{(n+2k-l-1)/2}\epsilon H^{-1}
                                        \beta_{n}(A) \eta_{n-1,k,l+1},$}\\
&K\eta_{nkl}=\makebox[0cm][l]{ 
$q^{(-n-2k+l)/2}\epsilon H\eta_{nkl}$ \quad  
on \ \ $\Hh=\displaystyle\mathop{\oplus}_{n,k,l=-\infty}^\infty \Hh_{nkl},
\ \,\Hh_{nkl}=\K.$}\\[-24pt]
\end{align*}
\begin{align*}
(I&.3)_{A,H,v,\epsilon}:\quad  &
a\eta_{nk}&=\alpha_n(A)\eta_{n-1,k},\quad  &
d\eta_{nk}&=\alpha_{n+1}(A)\eta_{n+1,k}, \quad \quad  \\
& & b\eta_{nk}&=q^{n+1}  A\eta_{n,k+1},\qquad &  
c\eta_{nk}&=q^n A\eta_{n,k-1},\\
&F\eta_{nk}=\makebox[0cm][l]{$ \lambda^{-1}
             q^{(-n+k-1)/2} \epsilon Hv\eta_{n,k-2}
+\lambda^{-1}q^{(n-k+1)/2}\epsilon H^{-1}\beta_{n+1}(A)\eta_{n+1,k-1},$} \\
&E\eta_{nk}=\makebox[0cm][l]{$  
-\lambda^{-1} q^{(-n+k+1)/2} \epsilon Hv^\ast \eta_{n,k+2}
  -\lambda^{-1}q^{(n-k-1)/2}\epsilon H^{-1} \beta_{n}(A) \eta_{n-1,k+1},$}\\
&K\eta_{nk}=\makebox[0cm][l]{ $ q^{(-n+k)/2}\epsilon H \eta_{nk}$    
\quad on \ \ $\Hh=
\displaystyle\mathop{\oplus}_{n,k=-\infty}^\infty \Hh_{nk},\ \,\Hh_{nk}=\K.$}
\end{align*}
Here, $A$, $B$, $H$ are commuting self-adjoint operators 
acting on a Hilbert space $\K$ 
such that $\sigma(A)\sqsubseteq (q,1]$, $\sigma (B)\sqsubseteq [q,1)$,
and $\sigma(H)\sqsubseteq (q^{1/2},1]$.
In the last series, 
$v$ is a unitary operator on $\K$ satisfying $Av=v A$ and $Hv=v H$. 
The parameter $\epsilon$ takes values in $\{-1,1\}$. 
Rep\-re\-sen\-ta\-tions 
labeled by different sets of parameters (within unitary equivalence)
or belonging to 
different series are not unitarily equivalent. 
A rep\-re\-sen\-ta\-tion of 
this series is irreducible if and only if $\K=\C$. In this case, we 
can regard the parameters $A$, $H$, $B$ and $v$ as complex numbers such that 
$A\in (q,1]$, $H\in(q^{1/2},1]$, $B\in[q,1)$ and $|v|=1$.
%
%
\subsection{Representations of the {\mathversion{bold}$\ast$}-al\-ge\-bra 
{\mathversion{bold}$\suu\!\ltimes\!\qd$} }
%
                                                    \label{rep-qd}
Clearly, $\qd$ is a $\ast$-sub\-al\-ge\-bra of\, $\suu\lti\qd$. 
By Lemma \ref{L3}, there are three series of 
$\ast$-rep\-re\-sen\-ta\-tions of $\qd$. 
Our aim is to extend these  representations 
to $\ast$-rep\-re\-sen\-ta\-tions of the 
cross product algebra\, $\suu\lti\qd$. 

To begin, let us determine the action of $K$ 
on the Hilbert space $\Hh$ from Lemma \ref{L3}. 
Assuming that the commutation relations of $K$ with $E$, $z$  
and $z^\ast$ hold in strong sense, $\Hh$ 
decomposes into two reducing subspaces on which 
$K>0$ and $K<0$. 
Studying both cases separately, we can write 
$K=\epsilon |K|$, where $\epsilon\in \{1,-1\}$.
 
Note that $z^\ast$ and $K$ satisfy the same relations as $Q_0$ and 
$H_0$ in Equations \rf[T2] and \rf[T3]. If $z$ is given by the 
formulas of the series $(I)$ or $(II)_A$, then the same reasoning  
as in Subsection \ref{sec-SUU} shows that $K$ acts on 
$\Hh$ by $K\eta_n=q^n\epsilon H\eta_n$, where $H$ denotes an 
invertible positive self-adjoint operator on $\Hh_0$. 
In addition, we have $AH=HA$ in the case $(II)_A$. 

In the third series $(III)_v$, $z=v$ is a unitary operator. 
Thus we can apply Lemma \ref{1} to describe the representations 
of the relation $zK=q^{-1}zK$. It states that 
$\Hh=\oplus^\infty_{n=-\infty} \Hh_n$, 
where each $\Hh_n$ is $\Hh_0$, 
$z\eta_n=\eta_{n+1}$ and $K\eta_n=q^n\epsilon H\eta_n$, 
where $H$ denotes a self-adjoint operator on $\Hh_0$ such that 
$\sigma (H) \sqsubseteq (q,1]$.

In the first two series, it follows from $Sz^\ast z=z^\ast zS$,  
$Sz=zS$ and $SK=q KS$ by standard arguments used repeatedly in this paper 
that $S$ acts on $\Hh$ by $S\eta_n=S_0\eta_n$, where $S_0$ is an operator 
on $\Hh_0$ satisfying 
$$
S_0 H=qH S_0,\quad  S_0^\ast S_0-q^2 S_0S_0^\ast=1-q^2 
$$  
and, in the second series, $S_0A=AS_0$. 
Hence $S_0$, $S^\ast_0$ and $H^{-1}$ 
satisfy on $\Hh_0$ the same relations as $z$, $z^\ast$ and $K$ on $\Hh$ 
so that above results concerning $z$, $z^\ast$ and $K$ apply. 
The operator $A$ can be handled in the same way as the operator $A$ in 
Equations \rf[T2] and \rf[T3], where $Q_0^\ast$ plays the role of $S_0$. 
This determines the first two series of 
representations of $\suu\lti\qd$.

Consider now the third series. From $SK=qKS$, it follows 
that $S$ maps $\Hh_n$ into $\Hh_{n-1}$ since the 
relation is assumed to hold in strong sense.
Write $S\eta_n=S_n\eta_{n-1}$. The identity $Sz=zS$
implies $S_{n+1}=S_n$, therefore $S_n=S_0$ for all $n$. 
Applying $SK=qKS$ to vectors $\eta_n\in\Hh_n$ shows 
that $S_0H=HS_0$. 
On $\Hh_0$, we have again $S_0^\ast S_0-q^2S_0 S_0^\ast =1-q^2$. 
The representations of this relation are described in 
Lemma \ref{L3}. The operator $H$ is treated just as the operator $A$ 
in the preceding paragraph. This completes the discussion of the 
third series. 

Carrying out all details, 
we obtain the following nine series of 
$\ast$-rep\-re\-sen\-ta\-tions 
of $\U_q(\mathrm{su}_{1,1})\lti \qd$. Let $\K$ be a 
Hilbert space. Suppose that $H_i$, $A_i$, $i\hsp=\hsp 1,2$, 
are self-adjoint operators 
acting on $\K$ such that $\ker H_1\hsp=\hsp\{0\}$, 
$\sigma(H_2)\hsp\sqsubseteq \hsp(q,1]$ 
and $\sigma(A_i)\hsp\sqsubseteq \hsp(q^2,1]$, $i\hsp=\hsp 1,2$, 
and suppose that $v$ is a 
unitary operator on $\K$. Assume that 
$A_1A_2\hsp=\hsp A_2A_1$, $H_iA_j\hsp=\hsp A_jH_i$, $i,j\hsp=\hsp 1,2$, and 
$H_2v\hsp=\hsp vH_2$. 
Let $\epsilon\hsp\in\hsp\{\pm1\}$. 
The rep\-re\-sen\-ta\-tions will be labeled 
by $(I.1)_{H_1}$, $(II.1)_{A_1,H_1},{\dots}$ etc. Define the Hilbert 
spaces $\Hh={\oplus^\infty_{n,k=0}}\Hh_{nk}$ in the case 
$(I.1)_{H_1}$; $\Hh={\oplus^\infty_{n=0}}
{\oplus^\infty_{k=-\infty}}\Hh_{nk}$ 
in the cases 
$(I.2)_{A_2,H_1}$ and $(I.3)_{H_1,\epsilon}$; 
$\Hh={\oplus^\infty_{n=-\infty}}{\oplus^\infty_{k=0}}\Hh_{nk}$ 
in the cases $(II.1)_{A_1,H_1}$
and $(III.1)_{H_2,\epsilon}$; 
$\Hh={\oplus^\infty_{n,k=-\infty}}\Hh_{nk}$ 
in the cases $(II.2)_{A_1,A_2,H_1}$, $(II.3)_{A_1,H_1,\epsilon}$ 
and $(III.2)_{A_2,H_2,\epsilon}$; 
and $\Hh={\oplus^\infty_{n=-\infty}}\Hh_n$ in the case 
$(III.3)_{v,H_2,\epsilon}$; where each $\Hh_{nk}$ and each 
$\Hh_n$ is equal to $\K$. 
The operators $z$ and 
$z^\ast$ act as follows. 
\begin{align*}
&(I.1)_{H_1},\ (I.2)_{A_2,H_1},\ (I.3)_{H_2,\epsilon}: \quad 
z\eta_{nk}=\lambda_{n+1}\eta_{n+1,k},\quad
z^\ast\eta_{nk}=\lambda_n\eta_{n-1,k}, \\
&(II.2)_{A_1,H_1},\ (II.2)_{A_1,A_2,H_1},\ (II.3)_{A_1,H_2,\epsilon}: \quad
z\eta_{nk}=\alpha_{n+1}(A_1)\eta_{n+1,k}, \\
&z^\ast\eta_{nk}=\alpha_n(A_1)\eta_{n-1,k},\\
&(III.1)_{H_2,\epsilon},\  (III.2)_{H_2,A_2,\epsilon}:\quad
z\eta_{nk}=\eta_{n+1,k}, \quad 
z^\ast\eta_{nk}=\eta_{n-1,k}, \\
&(III.3)_{H_2,v,\epsilon}: \quad 
z\eta_n=\eta_{n+1},\quad z^\ast\eta_n=\eta_{n-1}. 
\end{align*}
The operators $E,F$ and $K$ are given by
\begin{align*}
(I.1)&_{H_1}: \qquad \quad  K \eta_{nk}=q^{n-k} H_1 \eta_{nk},\\
&F \eta_{nk} = 
q^{n-k-1/2} \lambda^{-1}\lambda_{k+1} H_1 \eta_{n,k+1} 
+ q^{-(n-k)+1/2} \lambda^{-1}\lambda_n H_1^{-1} \eta_{n-1,k},\\
&E \eta_{nk} 
= \makebox[0cm][l]{$ -q^{n-k+1/2} \lambda^{-1}\lambda_k H_1 \eta_{n,k-1} 
- q^{-(n-k)-1/2} \lambda^{-1}\lambda_{n+1} H_1^{-1} \eta_{n+1,k},$} 
\phantom{\mbox{$-\lambda^{-1}q^{n-k+1/2} H_1\lambda_k (A_2) \eta_{n,k-1} 
-\lambda^{-1} q^{-(n-k)-1/2} \lambda_{n+1} H_1^{-1} \eta_{n+1,k},$}}
\\[-18pt]         
\end{align*}
\begin{align*}
(I.2)&_{A_2,H_1}: \qquad K \eta_{nk}=q^{n-k} H_1 \eta_{nk},\\
& F \eta_{nk} = 
\lambda^{-1} q^{n-k-1/2} H_1 \alpha_{k+1} (A_2) \eta_{n,k+1} 
+ \lambda^{-1}q^{-(n-k)+1/2} \lambda_n H_1^{-1} \eta_{n-1,k}, \\
& E \eta_{nk} 
= -\lambda^{-1}q^{n-k+1/2} H_1\alpha_k (A_2) \eta_{n,k-1} 
-\lambda^{-1} q^{-(n-k)-1/2} \lambda_{n+1} H_1^{-1} \eta_{n+1,k},  
\\[-18pt]           
\end{align*}
\begin{align*}
(I.3)&_{H_2,\epsilon}:\qquad  \ \,
K \eta_{nk}=q^{n-k} \epsilon H_2 \eta_{nk},\\
& F \eta_{nk} 
= \lambda^{-1}q^{n-k-1/2} \epsilon  H_2 \eta_{n,k+1} 
+ \lambda^{-1}q^{-(n-k)+1/2} \lambda_n\epsilon H_2^{-1} \eta_{n-1,k},\\
&  E \eta_{nk} 
=\makebox[0cm][l]{$   -\lambda^{-1} q^{n-k+1/2} \epsilon H_2 \eta_{n,k-1} 
-\lambda^{-1}q^{-(n-k)-1/2} \lambda_{n+1}\epsilon H_2^{-1} \eta_{n+1,k},$}
\phantom{\mbox{$-\lambda^{-1}q^{n-k+1/2} H_1\lambda_k (A_2) \eta_{n,k-1} 
-\lambda^{-1} q^{-(n-k)-1/2} \lambda_{n+1} H_1^{-1} \eta_{n+1,k},$}}  
\\[-18pt]
\end{align*}
\begin{align*}
(&II.\makebox[0cm][l]{$ 1)_{A_1,H_1}: $}
\phantom{\mbox{$ 2)_{A_1,A_2,H_1}: $}} \qquad 
K \eta_{nk}=q^{n-k} H_1 \eta_{nk},\\
&F \eta_{nk} 
= \lambda^{-1}q^{n-k-1/2} \lambda_{k+1} H_1 \eta_{n,k+1} 
+ \lambda^{-1}q^{-(n-k)+1/2} H_1^{-1} \alpha_n (A_1)\eta_{n-1,k},\\
&E \eta_{nk} 
=\makebox[0cm][l]{$ -\lambda^{-1}q^{n-k+1/2} \lambda_k H_1 \eta_{n,k-1} 
- \lambda^{-1}q^{-(n-k)-1/2} H_1^{-1}\alpha_{n+1} (A_1) \eta_{n+1,k}, $}
\phantom{\mbox{$- \lambda^{-1}q^{n-k+1/2} H_1\alpha_k (A_2) \eta_{n,k-1} 
- \lambda^{-1}q^{-(n-k)-1/2} H_1^{-1} \alpha_{n+1} (A_1) \eta_{n+1,k},$}}
\\[-18pt]
\end{align*}
\begin{align*}
(&II.2)_{A_1,A_2,H_1}: \qquad  K \eta_{nk}=q^{n-k} H_1 \eta_{nk},\\
& F \eta_{nk} 
= \lambda^{-1} q^{n-k-1/2} H_1 \alpha_{k+1} (A_2) \eta_{n,k+1} 
+ \lambda^{-1}q^{-(n-k)+1/2} H_1^{-1} \alpha_n (A_1)\eta_{n-1,k},\\
&E \eta_{nk} 
=- \lambda^{-1}q^{n-k+1/2} H_1\alpha_k (A_2) \eta_{n,k-1} 
- \lambda^{-1}q^{-(n-k)-1/2} H_1^{-1} \alpha_{n+1} (A_1) \eta_{n+1,k},
\\[-18pt]
\end{align*}
\begin{align*}
(&II.\makebox[0cm][l]{$ 3)_{A_1,H_2,\epsilon}: $}
\phantom{\mbox{$ 2)_{A_1,A_2,H_1}: $}} \qquad 
 K \eta_{nk}=q^{n-k}\epsilon H_2 \eta_{nk},\\
&F \eta_{nk} 
= \lambda^{-1}q^{n-k-1/2} \epsilon H_2 \eta_{n,k+1} 
+ \lambda^{-1}q^{-(n-k)+1/2} \epsilon H_2^{-1}\alpha_{n}(A_1)\eta_{n-1,k},\\
&E \eta_{nk} 
= \makebox[0cm][l]{$ -\lambda^{-1}q^{n-k+1/2} \epsilon H_2 \eta_{n,k-1} 
-\lambda^{-1}q^{-(n-k)-1/2}\epsilon H_2^{-1}\alpha_{n+1}(A_1)\eta_{n+1,k},$}
\phantom{\mbox{$- \lambda^{-1}q^{n-k+1/2} H_1\alpha_k (A_2) \eta_{n,k-1} 
- \lambda^{-1}q^{-(n-k)-1/2} H_1^{-1} \alpha_{n+1} (A_1) \eta_{n+1,k},$}}
\\[-18pt]
\end{align*}
\begin{align*}
(III&.\makebox[0cm][l]{$ 1)_{H_2,\epsilon}: $} 
\phantom{\mbox{$ 2)_{A_2,H_2,\epsilon}: $}} \qquad 
K \eta_{nk}=q^{n}\epsilon H_2 \eta_{nk}, \\
&F \eta_{nk} 
=\makebox[0cm][l]{$  
\lambda^{-1}q^{n-1/2} \lambda_{k+1}\epsilon H_2 \eta_{n-1,k+1} 
+\lambda^{-1}q^{-n+1/2} \epsilon H_2^{-1} \eta_{n-1,k}, $} 
\phantom{\mbox{$ \lambda^{-1}q^{n-1/2} \epsilon H_2 \alpha_{k+1}(A_2)
\eta_{n-1,k+1}+\lambda^{-1}q^{-n+1/2}\epsilon H_2^{-1}\eta_{n-1,k},$}}\\
&E \eta_{nk} 
= -\lambda^{-1}q^{n+1/2}\lambda_k \epsilon H_2 \eta_{n+1,k-1} 
-\lambda^{-1}q^{-n-1/2}\epsilon H_2^{-1} \eta_{n+1,k},
\\[-18pt]
\end{align*}
\begin{align*}
(III&.2)_{A_2,H_2,\epsilon}:  \qquad 
K \eta_{nk}=q^{n}\epsilon H_2 \eta_{nk}, \\
&F \eta_{nk} 
= \lambda^{-1}q^{n-1/2} \epsilon H_2 \alpha_{k+1}(A_2)\eta_{n-1,k+1} 
 +\lambda^{-1}q^{-n+1/2} \epsilon H_2^{-1} \eta_{n-1,k}, \\
&E \eta_{nk} 
= -\lambda^{-1}q^{n+1/2} \epsilon H_2 \alpha_k(A_2)\eta_{n+1,k-1} 
- \lambda^{-1}q^{-n-1/2} \epsilon H_2^{-1} \eta_{n+1,k},\\[-18pt]
\end{align*}
\begin{align*}
(III&.\makebox[0cm][l]{$ 3)_{v,H_2,\epsilon}: $} 
\phantom{\mbox{$ 2)_{A_2,H_2,\epsilon}: $}} \qquad 
K \eta_{n}=q^{n}\epsilon H_2 \eta_{n},\\
&F \eta_n 
=\makebox[0cm][l]{$  \lambda^{-1}q^{n-1/2}\epsilon H_2 v \eta_{n-1} 
+\lambda^{-1}q^{-n+1/2}\epsilon H_2^{-1} \eta_{n-1}, $} 
\phantom{\mbox{$ \lambda^{-1}q^{n-1/2} \epsilon H_2 \alpha_{k+1}(A_2)
\eta_{n-1,k+1}+\lambda^{-1}q^{-n+1/2}\epsilon H_2^{-1}\eta_{n-1,k},$}}\\
&E \eta_{n} 
= -\lambda^{-1}q^{n+1/2}\epsilon H_2 v^\ast \eta_{n+1} 
-\lambda^{-1}q^{-n-1/2}\epsilon H_2^{-1} \eta_{n+1},
\end{align*}

\noindent
Rep\-re\-sen\-ta\-tions 
labeled by different sets of parameters (within unitary equivalence)
or belonging to 
different series are not unitarily equivalent. 
A rep\-re\-sen\-ta\-tion of this list is irreducible if and only if $\K=\C$. 
In this case, the parameters $A_i$, $H_i$, $i=1,2$, become real numbers 
such that $H_1\ne 0$, $H_2\in(q,1]$, $A_i \in (q^2,1]$, $i=1,2$, 
and $v$ becomes a complex number of modulus 1. 
%
%
\subsection{Heisenberg representations of the cross product algebra 
{\mathversion{bold}$
\cU_q (\mathrm{su}_{1,1})\!\ltimes\!\cO(\mathrm{SU}_q(1,1))$}}
                                                         \label{sec-Heis}
We proceed in a similar manner as in Subsection \ref{3.5}. 
Let $\cF (\mathrm{SU}_q (1,1))$ denote the $\ast$-al\-ge\-bra 
generated by the algebra $\dC[u,v]$ of complex Laurent polynomials 
in commuting variables $u$, $v$ and the algebra $\cF(\dR^+)$ of 
locally bounded Borel functions 
on $\dR^+=(0,+\infty)$ with cross commutation relations  
and involution given by 
$$
u^n v^k f(r) = f(q^kr) u^n v^k,\quad 
(u^n v^k f(r))^\ast = \bar{f} (r) v^{-k} u^{-n}, 
$$ 
where $n,k\in \dZ$ and $f\in \cF(\dR^+)$. Define
\begin{align*}
&u^n v^k f(r) \anf E\hsp=\hsp
q^{\frac{n+k+1}{2}}\lambda^{-1} u^{n-1} v^{k+1} 
\Big(f(r) \sqrt{1 {+} q^{-2k} r^2} {-} q^{-n} f(q^{-1} r)
\sqrt{1 {+} r^2}\Big)r^{-1}\\
&u^n v^k f(r) \anf F\hsp =\hsp 
q^{\frac{k-n-3}{2}}\lambda^{{-}1} u^{n+1} v^{k-1} \Big( f(qr) 
\sqrt{1 {+} q^{2} r^2} {-} q^{n} f(r)\sqrt{1 {+} q^{-2k+2}r^2}\Big) r^{-1}\\
&u^n v^k f(r)\anf K \hsp=\hsp q^{\frac{n-k}{2}} u^n v^k f(r).
\end{align*}
Straightforward computations show 
that these formulas define indeed a right action of the 
Hopf $\ast$-al\-ge\-bra $\cU_q(\mathrm{su}_{1,1})$ 
on $\cF (\mathrm{SU}_q(1,1))$ 
such that the $\ast$-al\-ge\-bra $\cF (\mathrm{SU}_q(1,1))$ 
becomes a right $\cU_q (\mathrm{su}_{1,1})$-mod\-ule $\ast$-al\-ge\-bra. 
We omit the details of this lengthy and tedious  verification.

Further, one easily checks that there is an injective 
$\ast$-ho\-mo\-mor\-phism $\phi$ from $\cO (\mathrm{SU}_q(1,1))$ 
into $\cF (\mathrm{SU}_q(1,1))$ given by $\phi (a) = v \sqrt{1+r^2}$ 
and  $\phi (c) = ur$ such that $x\anf f = \phi (x) \anf f$ for 
$x \in \cO(\mathrm{SU}_q(1,1))$ and $f\in \cU_q (\mathrm{su}_{1,1})$. 
We shall identify $x\in \cO(\mathrm{SU}_q(1,1))$ with
$\phi (x) \in\cF (\mathrm{SU}_q(1,1))$. Then $\cO (\mathrm{SU}_q(1,1))$ 
is a right $\cU_q (\mathrm{su}_{1,1})$-mod\-ule $\ast$-sub\-al\-ge\-bra 
of $\cF (\mathrm{SU}_q(1,1))$ and the cross product 
algebra $\cU_q(\mathrm{su}_{1,1})\lti \cO (\mathrm{SU}_q(1,1))$ 
is a $\ast$-sub\-al\-ge\-bra of 
$\cU_q(\mathrm{su}_{1,1})\lti \cF (\mathrm{SU}_q(1,1))$.
Using the identities 
$a\hsp =\hsp v \sqrt{1+r^2}$,\ \,$d\hsp =\hsp a^\ast 
\hsp = \hsp \sqrt{1+r^2} u^{-1}$,\ \,$b\hsp  = \hsp 
q c^\ast\hsp  =\hsp  q u^{-1} r$,\ \,$c=ur$, 
it follows from the definition of the action $\anf$ on 
$\cF (\mathrm{SU}_q(1,1))$ that 
\begin{align*}
&a^k b^l c^n f(r) \anf E \hsp=\hsp
q^{\frac{k+n-l+1}{2}}
\lambda^{-1} a^{k+1} b^l c^{n-1} \big(f(r){-}q^{-2n} f(q^{-1} r) \big),\\
&d^k b^l c^n f(r) \anf E \hsp=\hsp
q^{\frac{n-k-l+1}{2}}\lambda^{-1} d^{k+1} b^l c^{n-1} 
\big((1 {+} q^{-2k}r^2)f(r) {-} (1 {+} r^2)q^{-2n}f(q^{-1} r)\big),
\end{align*}
for $k,l,n\in \dN_0$ and $f\in \cF(\dR^+)$. 
If we set $f\equiv 1$ and $r^2= q^{-1} bc$, we recover 
the action of $E$ on the vector space basis 
$\{ a^k b^l c^n,\  d^j b^l c^n\,;\,k,l,n \in \dN_0,\ j\in\dN\}$ 
of $\cO(\mathrm{SU}_q(1,1))$. 
The  action of $F$ on this basis is easily obtained 
from action of $E$ by applying 
$x\anf F= q^{-1} x\anf S(E)^\ast = q^{-1} (x^\ast \anf E)^\ast$. 

The construction  of a $\cU_q (\mathrm{su}_{1,1})$-in\-var\-i\-ant 
linear functional $h_{\mu_0}$ and of the corresponding 
Heisenberg representation is completely similar to Subsection \ref{3.5}. 
We fix a finite positive Borel measure $\mu_0$ on $(q,1]$ and extend it to 
a Borel measure $\mu$ on $\dR^+$ such that 
$\mu(q^k\cM)=q^k \mu_0 (\cM)$ for $k\in\dZ$ and $\cM\subseteq (q,1]$.
Let $\cF_0 (\mathrm{SU}_q(1,1))$ be the subalgebra of 
$\cF(\mathrm{SU}_q(1,1))$ generated by the elements 
$p(u,v)f(r)$, where $p(u,v)\in\dC[u,v]$ and $f(r)\in \cF(\dR^+)$
has compact support. 
Then the formula
$$
h_{\mu_0} (p(u,v) f(r)) 
= \into_{\dT^2} p(u,v)\dd u \dd v {\into^\infty_0} f(r)r\,\dd\mu (r),
$$
defines a $\cU_q(\mathrm{su}_{1,1})$-in\-var\-i\-ant positive 
linear functional $h_{\mu_0}$ on the 
right $\cU_q(\mathrm{su}_{1,1})$-mod\-ule 
$\ast$-al\-ge\-bra $\cF_0 (\mathrm{SU}_q(1,1))$. 

The Heisenberg representation $\pi_h$ of 
$\cU_q(\mathrm{su}_{1,1})\lti \cF (\mathrm{SU}_q(1,1))$
associated with $h\equiv h_{\mu_0}$ acts on the Hilbert space 
$\cL^2 (\dT^2)\otimes \cL^2 (\dR^+, r\dd\mu)$. 
The actions of the generators 
$a, b, c, d\hsp\in\hsp\cO (\mathrm{SU}_q(1,1))$, 
$f(r)\in\cF (\R^+)$ 
and $X\in\suu$ are given by 
\begin{align*}
&\pi_h(a)(u^n v^k \zeta (r))=u^n v^{k+1}\alpha_{-k}(r)\zeta (r),\ \   
\pi_h (b) (u^n v^k \zeta (r)) = q^{-k+1} u^{n-1} v^k r \zeta (r),\\
&\pi_h(d)(u^n v^k \zeta (r))=u^n v^{k-1}\alpha_{-k+1}(r)\zeta (r),\ \ 
\pi_h (c) (u^n v^k \zeta (r)) = q^{-k} u^{n+1} v^k r \zeta (r),\\
&\hspace{94pt}
\pi_h(f(r))(u^n v^k \zeta (r)) =  u^{n} v^k f(q^{-k}r) \zeta (r),\\
&\hspace{94pt}
\pi_h (X) (u^n v^k \zeta (r)) = (u^n v^k \zeta (r))\anf S^{-1}(X).
\end{align*}

Let $\Hh={\oplus^\infty_{n,k,l=-\infty}}\Hh_{nkl}$, where each 
Hilbert space $\Hh_{nkl}$ is  $\cL^2 ((q,1], r\dd\mu_0)$. 
Define a linear operator 
$W:\Hh\rightarrow \cL^2 (\dT^2)\otimes \cL^2 (\dR^+, r\dd\mu)$ by 
$$
   W\zeta_{nkl}:=q^k u^{-l}v^{-n}\zeta(q^kr),\quad 
\zeta\in \cL^2 ((q,1], r\dd\mu_0),\ \, n,k,l\in\Z. 
$$
The reasoning from Subsection \ref{3.5} which shows that 
$W$ is unitary applies verbatim. 
Hence the Heisenberg representation on 
$\cL^2 (\dT^2)\otimes \cL^2 (\dR^+, r\dd\mu)$ is unitarily equivalent to 
a $\ast$-re\-pre\-sen\-ta\-tion on $\Hh$ determined by the following 
formulas:
\begin{align*}
&\hspace{71pt} a\zeta_{nkl}=\alpha_{n-k}(Q)\zeta_{n-1,kl},\quad 
b\zeta_{nkl}=q^{n-k+1}  Q\zeta_{nk,l+1},\\
&\hspace{71pt}d\zeta_{nkl}=\alpha_{n-k+1}(Q)\zeta_{n+1,kl},\quad 
c\zeta_{nkl}=q^{n-k} Q\zeta_{nk,l-1},\\
&F\zeta_{nkl}=\lambda^{-1}\big( q^{(n-l+1)/2} \beta_{n-k+1}(Q)\zeta_{n+1,k,l-1}
-q^{(-n+l-3)/2}\alpha_{k}(Q^{-1})\zeta_{n+1,k+1,l-1}\big),\\
&E\zeta_{nkl}=\lambda^{-1}\big( 
q^{(-n+l-1)/2}\alpha_{k-1}(Q^{-1})\zeta_{n-1,k-1,l+1} 
-q^{(n-l-1)/2}\beta_{n-k}(Q)\zeta_{n-1,k,l+1}\big),\\
&K\zeta_{nkl}=q^{(-n+l)/2}\zeta _{nkl},
\end{align*}
where $Q$ is the multiplication operator on $\cL^2 ((q,1], r\dd\mu_0)$.
As a sample, we verify the formula for the action of $E$ on $\Hh$
and compute  
\begin{align*}
&E\zeta_{nkl}=W^{-1}(-q^{k-1}(u^{-l}v^{-n}\zeta(q^kr))\anf E)\\
&= \hspace{-1pt} W^{-1}\big(
q^{-(n+l+1)/2}\lambda^{-1} q^k u^{-l-1} v^{-n+1} 
( q^{l} \sqrt{1 {+} r^{-2}}\zeta(q^{k-1} r)  
\hspace{-0.5pt} -\hspace{-0.5pt} 
\sqrt{q^{2n} {+}r^{-2}  }\zeta(q^k r)  )\big)\\
&= \hspace{-1pt} \lambda^{-1}\big(
q^{(-n+l+1)/2} \alpha_{k-1}(Q^{-1})\zeta_{n-1,k-1,l+1}
 - q^{(n-l-1)/2}\beta_{n-k}(Q)\zeta_{n-1,k,l+1}\big). 
\end{align*}
Applying the unitary transformation 
$U\eta_{nkl}:=(-1)^k\zeta_{n+k+2,k+2,l-k+2}$ 
and rewriting above formulas in terms 
of $\eta_{nkl}$ proves the following proposition. 
\begin{thp}
The Heisenberg representation 
of $\cU_q(\mathrm{su}_{1,1})\lti \cO (\mathrm{SU}_q(1,1))$ 
associated with $h_{\mu_0}$ 
is unitarily equivalent 
to the representation $(II)_{Q,qQ^{-1},I,1}$ from 
Subsection \ref{sec-SUU}, where $Q$ 
denotes the multiplication operator on the Hilbert space 
$\cL^2 ((q,1], r\dd\mu_0)$. 
\end{thp}
%
%
\subsection{Heisenberg representations of the cross product algebra
{\mathversion{bold}$\cU_q(\mathrm{su}_{1,1})\!\ltimes\!\cO(\mathrm{U}_q)$}}
%
Analogously to Subsections \ref{3.5} and \ref{sec-Heis}, we denote by 
$\cF(\dR{\setminus}\{0\})$ the $\ast$-al\-ge\-bra of all 
locally bounded Borel functions on $\dR{\setminus}\{0\}$. 
Let $\cF(\mathrm{U}_q)$ be the $\ast$-al\-ge\-bra generated 
by the two $\ast$-al\-ge\-bras $\cO(\mathrm{U}_q)$ and 
$\cF(\dR{\setminus}\{0\})$ 
with relations
$$
z^\ast z=1 - y,\quad
z^nz^{\ast k} f(y)=f(q^{2(k-n)}y) z^n z^{\ast k},\quad 
n, k\in \dN_0,\ \, f\in \cF (\dR{\setminus}\{0\}).
$$
From the defining relations of $\cF(\mathrm{U}_q)$,  
it follows that each element $x\in\cF(\mathrm{U}_q)$ can be written as 
a finite sum 
\begin{equation}\label{Zeh}
x=\sumop_{n\ge 1} z^n f_n (y) + f_0 (y) 
+ \sumop_{k\ge 1} f_{-k} (y) z^{\ast k},
\end{equation}
where the functions $f_j(y)\in\cF(\dR{\setminus}\{0\})$ are uniquely 
determined by $x$.

We define a right action $\anf$ 
of the Hopf algebra $\cU_q(\mathrm{su}_{1,1})$ on $\cF(\mathrm{U}_q)$ by
\begin{align}\label{Nagel}
&z^n f(y) \anf E 
= q^{1/2} \lambda^{-1} z^{n+1} (q^{-n} f(y) - q^n f(q^2 y)),\\
\label{Ferse}
&f(y)z^{\ast k} \anf E 
= q^{1/2} \lambda^{-1} [zq^k (f(y) - f(q^2 y)) z^{\ast k} 
+ (q^k-q^{-k}) f(y)z^{\ast k-1}],\\
\label{Fessel}
&z^k f(y) \anf F 
= q^{-1/2}\lambda^{-1}  [z^k q^k (f(y) - f(q^2 y)) z^{\ast} 
+ (q^k-q^{-k}) z^{k-1} f(y)],\\
\label{Wade}
&f(y) z^{\ast n} \anf F = q^{-1/2} \lambda^{-1} (q^{-n} f(y) 
- q^n f(q^2 y))z^{\ast(n+1)},\\
\label{Knie}
&z^n f(y) \anf K = q^{-n} z^n f(y),\quad  
f(y) z^{\ast k} \anf K=q^{-k} f(y) z^{\ast k}
\end{align}
for $n\in \dN_0, k\in \dN$ and $f\in \cF(\dR{\setminus}\{0\})$. 
Then $\cF(\mathrm{U}_q)$ is 
a right $\cU_q(\mathrm{su}_{1,1})$-mod\-ule $\ast$-al\-ge\-bra 
and $\cO(\mathrm{U}_q)$ is 
a right $\cU_q(\mathrm{su}_{1,1})$-mod\-ule $\ast$-sub\-al\-ge\-bra 
of $\cF(\mathrm{U}_q)$. As in the previous subsections, 
we omit the details of these verifications.

Our next aim is to define a $\cU_q(\mathrm{su}_{1,1})$-in\-var\-i\-ant 
positive linear functional on an appropriate 
$\cU_q(\mathrm{su}_{1,1})$-mod\-ule $\ast$-sub\-al\-ge\-bra 
of $\cF(\mathrm{U}_q)$. 
Let $\cF_0(\mathrm{U}_q)$ denote the $\ast$-al\-ge\-bra generated 
by the elements $z^kz^{\ast n} f(y)$, 
where $k,n\in\dN_0$ and $f(y)\in \cF(\dR{\setminus}\{0\}) $ 
has  compact support.  
With the $\suu$-ac\-tion given by \rf[Nagel]--\rf[Knie], 
$\cF_0(\mathrm{U}_q)$ becomes a $\suu$-mod\-ule $\ast$-al\-ge\-bra. 
Let $h$ be a linear functional on $\cF_0(\mathrm{U}_q)$. 
From (\ref{Knie}), it follows
that $h$ is invariant under $K$ if and only if
\begin{equation}\label{Schenkel}
h(z^n f(y)) = h(f(y) z^{\ast n}) = 0\quad \mbox{for\ all}\ \, n\in\dN.
\end{equation}
If (\ref{Schenkel}) holds, 
then (\ref{Nagel})--(\ref{Wade}) imply 
that $h$ is invariant under $E$ and $F$ if 
and only if $h(f(y) z^\ast \anf E) = h(z f(y) \anf F) =0$. 
By (\ref{Ferse}) and (\ref{Fessel}), 
the latter is equivalent to 
\begin{equation}\label{Knochen}
h(f(q^{-2} y) (1 - q^{-2} y)) = q^{-2} h(f(y) (1 - y)).
\end{equation}
That is, $h$ is invariant if and only if (\ref{Schenkel}) 
and (\ref{Knochen}) are satisfied. Now we turn to the positivity condition. 
Let $x$ be as in (\ref{Zeh}). 
The relations $z^\ast z - q^2 zz^\ast 
= 1-q^2$ and $z^\ast z = 1 - y$ imply that
$$
z^{\ast n} z^n 
= \mathop{\mbox{$\prod$}}^{n-1}_{l=0} (1 - q^{2l} y),\quad z^nz^{\ast n} 
= \mathop{\mbox{$\prod$}}^{n}_{l=1} (1 - q^{-2l} y),\quad n\in\dN.
$$
Combining the latter with condition (\ref{Schenkel}), 
we obtain
\begin{multline}\label{Bauch}
h(x^\ast x) 
=\sumop_{n\ge 1} h
\Big( |f_n(y)|^2 \mathop{\mbox{$\prod$}}^{n-1}_{l=0} (1 - q^{2l} y)\Big) 
 +  h(|f_0 (y)|^2)   \\
+\sumop_{k\ge 1} h
\Big( |f_{-k} (y)|^2 \mathop{\mbox{$\prod$}}^{k}_{l=1}  (1 - q^{-2l} y)\Big).
\end{multline}
Let us suppose that $h$ is given by 
a positive measure on $\dR{\setminus}\{0\}$. 
From (\ref{Bauch}), we conclude that $h(x^\ast x)\ge 0$ for all $x$ 
provided that the support of the measure is contained 
in the set $(-\infty,0)\cup \{q^{2n}; n\in\dN\}$. 
We write the measure as a sum of measures 
with supports contained in $(-\infty,0)$ and $\{q^{2n};n\in\dN\}$, 
respectively. That is, we write $h=h_I+h_{\mu_0}$, 
where $h_I$ and $h_{\mu_0}$ are defined by 
\begin{align*}
&h_I (z^n f(y)) = h_I (f(y) z^{\ast n})
=\delta_{n0} \sumop^\infty_{k=1} f(q^{2k}) q^{-2k},\\
&h_{\mu_0} (z^n f(y)) = h_{\mu_0} (f(y) z^{\ast n})
=\delta_{n0} \into^0_{-\infty} f(y) y^{-2}\, \dd\mu(y).
\end{align*}
Here $\mu_0$ is a {\it finite} positive 
Borel measure on $[-1,-q^2)$ and $\mu$ denotes 
its extension to a Borel measure on $(-\infty,0)$ 
such that 
$\mu(q^{2k}\cM)=q^{2k} \mu_0(\cM)$ for $k\in\dZ$, 
$\cM\subseteq [-1,-q^2)$. 
One easily checks that 
Equation (\ref{Knochen}) is satisfied 
for both functionals $h_I$ and $h_{\mu_0}$. 
Hence $h_I$ and $h_{\mu_0}$ are $\cU_q(\mathrm{su}_{1,1})$-in\-var\-i\-ant 
positive linear functionals 
on the right $\cU_q(\mathrm{su}_{1,1})$-mod\-ule $\ast$-al\-ge\-bra 
$\cF_0 (\mathrm{U}_q)$.

Let $h=h_I$ or $h=h_{\mu_0}$. 
For the generators $z$ and $z^\ast$,  
the Heisenberg representation $\pi_h$ acts 
on the domain $\cF_0 (\mathrm{U}_q)$ by
\begin{align*}
&\pi_h (z) x = \sumop_{n\ge 0} z^{n+1} f_n (y) 
+ \sumop_{k\ge 1} f_{-k} (q^{-2} y) (1 - q^{-2} y) z^{\ast (k-1)},\\
&\pi_h (z^\ast) x = \sumop_{n\ge 1} z^{n-1}(1 - q^{2n-2}y) f_n (y) 
+ \sumop_{k\ge 0} f_{-k} (q^{2} y) z^{\ast (k+1)},
\end{align*}
where $x\in \cF_0(\mathrm{U}_q)$ is given by (\ref{Zeh}).

Let $u$ be the partial isometry 
from the polar decomposition of the closure 
of the operator $\pi_h(z)$. When $h=h_{\mu_0}$,  
the operator $u$ is a bilateral shift, so $u$ is unitary. 
Using this fact, we present another approach 
to the Heisenberg representation $\pi_h$. 
For notational convenience, we replace  $y$  by $-y$.

Let $\cF^+ (\mathrm{U}_q)$ be the $\ast$-al\-ge\-bra 
generated by the algebra $\dC[u]$ of complex Laurent polynomials in $u$ 
and the algebra $\cF (\dR^+)$ of locally bounded Borel 
functions  on $\dR^+$ with cross 
commutation relation  $u^n f(y) = f(q^{-2n} y) u^n$ 
and involution $(u^n f(y))^\ast = \bar{f} (y) u^{-n}$,  
where $n\in\dZ$ and $f\in\cF (\dR^+)$. 
Then $\cF^+ (\mathrm{U}_q)$ is 
a right $\cU_q(\mathrm{su}_{1,1})$-mod\-ule $\ast$-al\-ge\-bra 
with respect to the action $\anf$ given by the formulas
\begin{align*}
&u^n f(y) \anf E 
= q^{1/2} \lambda^{-1} u^{n+1} (q^{-n} f(y) \sqrt{1+q^{2n} y} 
- q^n f(q^2 y) \sqrt{1+y} ),\\
&u^n f(y) \anf F 
= q^{-1/2} \lambda^{-1} u^{n-1} (q^{n} f(q^{-2} y) \sqrt{1+q^{-2} y} 
- q^{-n} f(y) \sqrt{1+q^{2n-2} y} ),\\
&u^n f(y) \anf K = q^{-n} u^{n} f(y) 
\end{align*}
for $n\in\dZ$ and $f\in\cF(\dR^+)$. 
There is an injective $\ast$-ho\-mo\-mor\-phism 
$\phi : \cO(\mathrm{U}_q)\rightarrow 
\cF^+ (\mathrm{U}_q)$ such that $\phi (z) = u\sqrt{1+y}$ 
and $x\anf f = \phi (x)\anf f$ for $x\in \cO(\mathrm{U}_q)$ 
and $f\in \cU_q (\mathrm{su}_{1,1})$.
Thus we can consider $\cO(\mathrm{U}_q)$ 
as a $\cU_q(\mathrm{su}_{1,1})$-mod\-ule $\ast$-sub\-al\-ge\-bra 
of $\cF^+ (\mathrm{U}_q)$ by identifying $x\in\cO(\mathrm{U}_q)$ 
with $\phi(x)$.

Let $\mu_0$ be a finite positive Borel measure on $(q^2,1]$. 
We extend $\mu_0$ to a Borel measure $\mu$ on $\dR^+$ such 
that $\mu(q^{2k}\cM)= q^{2k} \mu_0 (\cM)$ for $k\in\dZ$, 
$\cM\subseteq (q^2,1]$. 
Let $\cF^+_0(\mathrm{U}_q)$ denote the $\ast$-sub\-al\-ge\-bra 
of $\cF^+ (\mathrm{U}_q)$ generated by the elements $u^k f(y)$, 
where $k\in\dZ$ and $f\in\cF(\dR^+)$ has  compact support. 
Define a linear functional $h_{\mu_0}$ on $\cF^+_0 (\mathrm{U}_q)$ by 
$$
\hat{h}_{\mu_0}(p(u) f(y))
=\into_\dT p(u) \dd u \into^{+\infty}_0 f(y) y^{-2}\,\dd\mu(y).
$$
From the definitions of the actions of the generators $E$, $F$, $K$, 
it follows immediately 
that $\hat{h}_{\mu_0} (x\anf f)=\varepsilon (f) \hat{h}_{\mu_0} (x)$ 
for $x\in \cF^+_0 (\mathrm{U}_q)$ and $f= E,F,K$. 
If $x=\sum_k u^k f_k(y)$, then
\begin{multline*}
\hat{h}_{\mu_0} (x^\ast x) 
= \hat{h}_{\mu_0} \Big(\sumop_{k,l} 
\overline{f}_k(y) u^{-k} u^l f_l (y)\Big)\\
=\sumop_{k,l} \hat{h}_{\mu_0} (u^{l-k} 
\overline{f_k} (q^{2(l-k)}y) f_l (y))
= \sumop_k \into^{+\infty}_0 |f_k(y)|^2 y^{-2}\, \dd\mu(y)\ge 0.
\end{multline*}
Therefore, $\hat{h}_{\mu_0}$ is a $\cU_q(\mathrm{su}_{1,1})$-in\-var\-i\-ant 
positive linear functional on $\cF^+_0 (\mathrm{U}_q)$.
If $\hat{\mu}_0 (\cM) = \mu_0 (-\cM)$ for all $\cM\subseteq (q^2,1]$,
then the Heisenberg representations associated with  $\hat{h}_{\mu_0}$ 
and $h_{\hat{\mu}_0}$ are unitarily equivalent.

Since $h=h_I+h_{\mu_0}$ is given by a direct sum of measures 
with supports contained in $(-\infty,0)$ and $\{q^{2n};n\in\dN\}$, and since 
the action of the cross product algebra 
$\cU_q(\mathrm{su}_{1,1})\lti \cF(\mathrm{U}_q)$
respects this decomposition, the Heisenberg 
representation associated with $h$ decomposes into a direct sum 
of Heisenberg representations associated with $h_I$ and $h_{\mu_0}$. 

Recall 
that the Heisenberg representations 
of $\cU_q(\mathrm{su}_{1,1})\lti \cF(\mathrm{U}_q)$
associated with $h_I$ 
acts on $\tilde{\cF}_0(\mathrm{U}_q):=\cF_0(\mathrm{U}_q)/\cN_I$, 
where $\cN_I:=\{ x\hsp\in\hsp\cF_0(\mathrm{U}_q)
\hspace{1.5pt};\hspace{1.5pt} h_I(x^\ast x)\hs{=}\hs 0\}$. 
For $k\in\N_0$, let $\delta_k(r)$ denote the characteristic function of 
the point $\{q^{2k}\}$, that is, $\delta_k(r)=1$ if $r=q^{2k}$ and 
zero otherwise. 
Note that $z^{\ast n}\delta_k(y)=\delta_k(q^{2n}y)z^{\ast n}=0$ 
in $\tilde{\cF}_0(\mathrm{U}_q)$ for $n\geq k$.  
Set 
\begin{align*}
&\zeta_{nk}:=
\big(\mathop{\mbox{$\prod$}}_{l=0}^{n-1}(1-q^{2(k+l)})\big)^{-1/2} 
q^k z^n \delta_k(y),\quad \zeta_{0k}:= q^k \delta_k(y),\quad
k\in\N, \  \,n\in\N,\\
&\zeta_{nk}:=
\big(\mathop{\mbox{$\prod$}}_{l=1}^{|n|}(1-q^{2(k-l)})\big)^{-1/2} 
q^k z^{\ast n}\delta_k(y),\quad k>1, \ \, n=-k+1,\ldots,-1.
\end{align*}
Then $\{\zeta_{nk}\,;\, k\in\N,\ n=-k+1,-k+2,\ldots\}$ is a set of 
orthonormal vectors which span $\tilde{\cF}_0(\mathrm{U}_q)$. 
Computing the actions of $z$, $z^\ast$,  $E$, $F$, $K$
and $f(r)\in \cF(\R{\setminus}\{0\})$ on $\zeta_{nk}$ gives 
\begin{align*}
&z\zeta_{nk}=\lambda_{n+k}\zeta_{n+1,k}, \quad 
z^\ast\zeta_{nk}=\lambda_{n+k-1}\zeta_{n-1,k}, \quad
f(r)\zeta_{nk}=f(q^{2(n+k)})\zeta_{nk}\\
&E\zeta_{nk}=\lambda^{-1}q^{n+1/2}\lambda_{k-1}\zeta_{n+1,k-1}
-\lambda^{-1}q^{-n-1/2}\lambda_{n+k}\zeta_{n+1,k},\\
&F\zeta_{nk}=-\lambda^{-1}q^{n-1/2}\lambda_{k}\zeta_{n-1,k+1}
+\lambda^{-1}q^{-n+1/2}\lambda_{n+k-1}\zeta_{n-1,k}, \quad 
K\zeta_{nk}= q^{n}\zeta_{nk}. 
\end{align*}
For instance, 
since $\delta_k(q^{2n}y)=\delta_{k-n}(y)$ 
and $z\delta_k(y)z^\ast=(1-q^{-2}y)\delta_{k+1}(y)
=(1-q^{2k})\delta_{k+1}(y)$, we have 
\begin{align*}
&z^{\ast m}\delta_k(y)\anf E=\delta_{k-m}(y)z^{\ast m}\anf E
=q^{1/2}\lambda^{-1}[q^m(1-q^{2(k-m)})\delta_{k-m+1}(y)\\
&\hspace{100pt}-q^m(1-q^{2(k-m-1)})\delta_{k-m}(y)
+(q^m-q^{-m})\delta_{k-m}(y)]z^{\ast m-1}\\
&=q^{1/2}\lambda^{-1}[q^m(1-q^{2(k-m)})z^{\ast m-1}\delta_{k}(y)
-q^{-m}(1-q^{2(k-1)})z^{\ast m-1}\delta_{k-1}(y)]. 
\end{align*}
From the latter expression, we derive the action of $E$ on $\zeta_{-m,k}$, 
$k,m\in\N$, $m<k$, by inserting the definition of 
the vectors $\zeta_{nk}$. 
If we set $\eta_{nk}:=(-1)^k\zeta_{n-k,k+1}$, $n,k\in\N_0$, 
then the actions of $z$, $z^\ast$, $K$, $E$ and $F$ on $\eta_{nk}$ 
coincide with the formulas of the series $(I.1)_{I}$ from 
Subsection \ref{rep-qd}.

Next we turn to the 
the Heisenberg representation 
of $\cU_q(\mathrm{su}_{1,1})\lti \cF(\mathrm{U}_q)$
associated with $h_{\mu_0}$. 
As noted above, it is 
unitarily equivalent to the Heisenberg representation 
associated with $\hat h_{\hat \mu_0}$, where $\hat\mu_0$ is a Borel 
measure on $(q^2,1]$ such that $\hat\mu_0(\cM)=\mu_0(-\cM)$ for 
$\cM\subset(q^2,1]$. The latter representation acts on the Hilbert space 
$\cL^2 (\dT)\otimes \cL^2 (\dR^+, y^{-2}\dd\mu)$ by 
\begin{align*}
&z(u^n \zeta(y))=u^{n+1}\sqrt{1+q^{2n}y}\hspace{2pt}\zeta(y),\quad 
z^\ast(u^n \zeta(y))=u^{n-1}\sqrt{1+q^{2(n-1)}y}\hspace{2pt}\zeta(y),\\
&f(y)(u^n \zeta(y))=u^n f(q^{2n}y)\zeta(y), \quad 
Z(u^n \zeta(y))=(u^n \zeta(y))\anf S^{-1}(Z),
\end{align*}
where $f(y)\in\cF(\R{\setminus}\{0\})$ and $Z\in\suu$. 

Let
$\Hh={\oplus^\infty_{n,k=-\infty}}\Hh_{nk}$, where each $\Hh_{nk}$
is  $\cL^2 ((q^2,1], y^{-2}\dd\mu_0)$, and consider the  linear operator 
$W:\Hh\rightarrow \cL^2 (\dT)\otimes \cL^2 (\dR^+, y^{-2}\dd\mu)$ defined by 
$$
   W\zeta_{nk}:=q^{k} u^{n}\zeta(q^{-2k}y),\quad 
\zeta\in \cL^2 ((q^2,1], y^{-2}\dd\mu_0),\ \, n,k,l\in\Z. 
$$
Similarly to Subsection \ref{sec-Heis}, one shows that $W$ is a unitary 
operator. Hence the Heisenberg representation associated with 
$\hat h_{\hat \mu_0}$ is unitarily equivalent to a 
$\ast$-re\-pre\-sen\-ta\-tion on $\Hh$. Let $Q$ denote 
the multiplication operator on $\cL^2 ((q^2,1], y^{-2}\dd\mu_0)$. 
Computing the actions of $z$, $z^\ast$, $K$, $E$ and $F$ on vectors 
$\zeta_{nk}$ gives 
\begin{align*}
&z\zeta_{nk}=\sqrt{1\hspace{-0.6pt}+\hspace{-0.6pt}q^{2(n+k)}Q}
\hs\zeta_{n+1,k},\ \,
z^\ast\zeta_{nk}=\sqrt{1\hspace{-0.6pt}+\hspace{-0.6pt}q^{2(n+k-1)}Q}
\hs \zeta_{n-1,k},\ \,
K\zeta_{nk}=q^n\zeta_{nk},\\
&E\zeta_{nk}=\lambda^{-1} q^{n+1/2}\sqrt{1+q^{2(k-1)}Q}\hs\zeta_{n+1,k-1}
-\lambda^{-1} q^{-n-1/2}\sqrt{1+q^{2(n+k)}Q\hs}\zeta_{n+1,k},\\
&F\zeta_{nk}=-\lambda^{-1} q^{n-1/2}\sqrt{1+q^{2k}Q}\hs\zeta_{n-1,k+1}
+\lambda^{-1} q^{-n+1/2}\sqrt{1+q^{2(n+k-1)}Q}\hs\zeta_{n-1,k}.
\end{align*}
Renaming $\eta_{nk}:=(-1)^k\zeta_{n-k,k+1}$ and computing the actions 
of $z$, $z^\ast$, $K$, $E$ and $F$ on $\eta_{nk}$, 
we obtain the formulas of the series $(II.2)_{Q,Q,1}$.

We summarize the preceding results in the following proposition. 
\begin{thp} 
The Heisenberg representation of 
$\cU_q(\mathrm{su}_{1,1})\lti \cO(\mathrm{U}_q)$
associated with $h$ 
is unitarily equivalent to the direct sum of the 
irreducible $\ast$-re\-pre\-sen\-ta\-tion $(I.1)_{I}$  
and the representation $(II.2)_{Q,Q,1}$ 
from Subsection \ref{rep-qd}, where $Q$ 
denotes the multiplication operator on 
$\K=\cL^2 ((q^2,1], y^{-2}\dd\mu_0)$. 
\end{thp}

\end{document}